\edef\savecatcodeat{\the\catcode`@}
\catcode`\@=11

\def\tb@ifSpecChars#1#2{#1}
\def\tb@ifNoSpecChars#1#2{#2}

\def\tableau{%
  \bgroup% matched in \tb@tableauD
  \@ifstar{\let\Tif\tb@ifNoSpecChars\tb@tableauB}% *, don't use special chars
          {\let\Tif\tb@ifSpecChars\tb@tableauB}}% no *, use special chars

\def\tb@tableauB{% add [] if no [options]
  \@ifnextchar[{\tb@tableauC}{\tb@tableauC[]}}

\def\tb@tableauC[#1]{\hbox\bgroup%
    \let\\=\cr% end line
    \def\bl{\global\let\tbcellF\tb@cellNF}%
    \def\tf{\global\let\tbcellF\tb@cellH}% highlighted cell
    \dimen2=\ht\strutbox \advance\dimen2 by\dp\strutbox%
    \ifx\baselinestretch\undefined\relax%
    \else%
       \dimen0=100sp \dimen0=\baselinestretch\dimen0%
       \dimen2=100\dimen2 \divide\dimen2 by\dimen0%
    \fi%
    \let\tpos\tb@vcenter% default position
    \tb@initYoung% default tableau type
    \tb@options#1\eoo% parse options
    \let\arrow\tb@arrow%
    \dimen0=\Tscale\dimen2%
    \dimen1=\dimen0 \advance\dimen1 by \tb@fframe%
    \lineskip=0pt\baselineskip=0pt% line spacing will be from \vbox to \dimen0
    \def\tb@nothing{}%
    \def\endcellno{$\rss\egroup\bss\egroup}% end cell w/o overlap
    \def\endcell{\endcellno\kern-\dimen0}% end cell & prepare to overlap it
    \def\begincell{\vbox to\dimen0\bgroup\vss\hbox to\dimen0\bgroup\hss$}%
    \let\overlay\tb@overlay%
    \let\fl\tb@fl%
    \let\lss\hss\let\rss\hss\let\tss\vss\let\bss\vss% cell alignment
    \def\mkcell##1{% format individual cell
        \let\tbcellF\tb@cellD% default cell frame
        \def\tb@cellarg{##1}% store cell contents
        \ifx\tb@cellarg\tb@nothing\let\tb@cellarg\tb@cellE\fi%
        \begincell\tb@cellarg\endcellno% the actual cell content
        \tbcellF}% draw cell frame
    \let\savecellF\tbcellF% save global value of cellF in case of nested tableau
     \Tif{\catcode`,=4\catcode`|=\active}{}\tb@tableauD}%

\let\tb@savetableauD\tableauD% save any current definition
{% set up characters which will be interpreted as command characters
    \catcode`|=\active \catcode`*=\active \catcode`~=\active%
    \catcode`@=\active% command characters
\gdef\tableauD#1{%
  \Tif{% make all the command characters active in math mode when #1 parsed
    \mathcode`|="8000 \mathcode`*="8000%
    \mathcode`~="8000 \mathcode`@="8000%
    \def@{\bullet}%
    \let|\cr% end line
    \let*\tf% highlighted cell
    \let~\sk% skew cell
  }{}%
  \tpos{\tabskip=0pt\halign{&\mkcell{##}\cr#1\crcr}}%
  \global\let\tbcellF\savecellF% restore global value
  \egroup% match \hbox\bgroup at start of \tableauC
  \egroup}% match \bgroup at start of \tableau
}
\let\tb@tableauD\tableauD% rename the command
\let\tableauD\tb@savetableauD% restore old command with this name
\let\tb@savetableauD\undefined

\def\tb@options#1{\ifx#1\eoo\relax\else\tb@option#1\expandafter\tb@options\fi}

\def\tb@option#1{%
  \if#1t\let\tpos\tb@vtop\fi%        t = align at top
  \if#1c\let\tpos\tb@vcenter\fi%     c = align at center
  \if#1b\let\tpos\vbox\fi%           b = align at bottom
  \if#1F\tb@initFerrers\fi%          F = Ferrers diagram
  \if#1Y\tb@initYoung\fi%            Y = Young diagram
  \if#1s\tb@initSmall\fi%            s = small boxes
  \if#1m\tb@initMedium\fi%           m = medium boxes
  \if#1l\tb@initLarge\fi%            l = large boxes
  \if#1p\tb@initPartition\fi%            p = small partition sized boxes
  \if#1a\tb@initArrow\fi%            a = use arrow font as base dimension
}

\def\tb@vcenter#1{\ifmmode\vcenter{#1}\else$\vcenter{#1}$\fi}

\def\tb@vtop#1{\hbox{\raise\ht\strutbox\hbox{\lower\dimen0\vtop{#1}}}}

\def\tb@initPartition{\def\Tscale{.3}}
\def\tb@initSmall{\def\Tscale{1}}
\def\tb@initMedium{\def\Tscale{2}}
\def\tb@initLarge{\def\Tscale{3}}

\def\tb@initArrow{\dimen2=1.25em}

\def\tb@initYoung{%
  \def\tb@cellE{}% empty cell stays empty
  \let\tb@cellD\tb@cellN% default frame is normal frame
  \def\sk{\global\let\tbcellF\tb@cellNF}}% skew cells are empty
\def\tb@initFerrers{%
  \def\tb@cellE{\bullet}% empty cell gets bullet
  \let\tb@cellD\tb@cellNF% default frame is no frame
  \def\sk{\bullet}}% skew cell gets bullet

\tb@initMedium% default scale

\def\tb@sframe#1{%
  \vbox to0pt{%            Embed frame in a box of no vert or hor extent
    \vss%                            pull box above reference point
    \hbox to0pt{%
      \hss%                          pull box left of reference point
      \vbox to\dimen1{%              Actual width of frame
        \hrule depth #1 height0pt% draw top edge of frame
        \vss%                     begin vcenter sides
        \hbox to\dimen1{%           horiz box with side edges just inside
          \vrule width #1 height\dimen1% left edge
          \hss%                     stretch center
          \vrule width #1%         right edge
          }%
        \vss%                     end vcenter sides
        \hrule height #1 depth 0in% bottom edge
        }%
      \kern-\tb@hframe%           horiz alignment off by half line width
      }%
    \kern-\tb@hframe}}%           vert alignment off by half line width
\def\tb@hframe{.2pt}\def\tb@fframe{.4pt}\def\tb@bframe{2pt}
\def\tb@cellH{\tb@sframe{\tb@bframe}}       % bold frame
\def\tb@cellNF{}                            % no frame
\def\tb@cellN{\tb@sframe{\tb@fframe}}       % normal frame
\let\tbcellF\tb@cellN                       % default is normal

\def\tb@rpad{1pt}
\def\tb@lpad{1pt}
\def\tb@tpad{1.8pt}
\def\tb@bpad{1.8pt}

\def\tb@overlay{\endcell\@ifnextchar[{\tb@overlaya}{\begincell}}
\def\tb@overlaya[#1]{\vbox to\dimen0\bgroup%
  \tb@overlayoptions#1\eoo%
  \tss\hbox to\dimen0\bgroup\lss$}
\def\tb@overlayoptions#1{\ifx#1\eoo\relax\else\tb@overlayoption#1\expandafter\tb@overlayoptions\fi}

\def\tb@overlayoption#1{
  \if#1t\def\tss{\vskip\tb@tpad}\let\bss\vss\fi% t = align at top
  \if#1c\let\tss\vss\let\bss\vss\fi%             c = align at center
  \if#1b\def\bss{\vskip\tb@bpad}\let\tss\vss\fi% b = align at bottom
  \if#1l\def\lss{\hskip\tb@lpad}\let\rss\hss\fi% l = align at left
  \if#1m\let\lss\hss\let\rss\hss\fi%             m = align at middle
  \if#1r\def\rss{\hskip\tb@rpad}\let\lss\hss\fi% r = align at right
}

\def\tb@fl{\endcell\begincell\vrule depth 0pt width \dimen0 height \dimen0 \endcell\begincell}

\@ifundefined{diagram}{}{
\def\tb@arrowpad{.5}

\newoptcommand{\tb@arrow}{\@ne}[2]{%
  \endcell% end previous cell contents
   \begingroup%
   \let\dg@getnodesize\tb@getnodesize% substitute routine to get nodesize
   \dg@USERSIZE=#1\relax%
   \ifnum\dg@USERSIZE<\@ne \dg@USERSIZE=\@ne \fi%
   \dg@parse{#2}%
   \dg@label{\tb@draw{#1}{#2}}}% draw arrow

\def\tb@getnodesize#1#2#3#4#5{\dimen3=\tb@arrowpad\dimen2 #4=\dimen3 #5=\dimen3\relax}
\def\tb@getnodesize#1#2#3#4#5{\ifnum#2=0\ifnum#3=0\tb@getnodesizetail{#4}{#5}\else\tb@getnodesizehead{#4}{#5}\fi\else\tb@getnodesizehead{#4}{#5}\fi}
\def\tb@getnodesizetail#1#2{\dimen3=.5\dimen2 #1=\dimen3 #2=\dimen3}
\def\tb@getnodesizehead#1#2{\dimen3=.5\dimen2 #1=\dimen3 #2=\dimen3}

\def\tb@draw#1#2#3#4{%
        \dg@X=0\dg@Y=0\dg@XGRID=1\dg@YGRID=1\unitlength=.001\dimen0%
        \dg@LBLOFF=\dgLABELOFFSET \divide\dg@LBLOFF\unitlength%
        \dg@drawcalc% compute arrow geometry
        \begincell% start tableau cell
        \let\lams@arrow\tb@lams@arrow% substitute routine
        \begin{picture}(0,0)\begingroup\dg@draw{#1}{#2}{#3}{#4}\end{picture}%
        \endcell% end tableau cell
        \endgroup% match \begingroup in \tb@arrow
        \begincell}% start new entry in this cell
}

\def\tb@lams@arrow#1#2{%
 \lams@firstx\z@\lams@firsty\z@
 \lams@lastx#1\relax\lams@lasty#2\relax
 \lams@center\z@
 \N@false\E@false\H@false\V@false
 \ifdim\lams@lastx>\z@\E@true\fi
 \ifdim\lams@lastx=\z@\V@true\fi
 \ifdim\lams@lasty>\z@\N@true\fi
 \ifdim\lams@lasty=\z@\H@true\fi
 \NESW@false
 \ifN@\ifE@\NESW@true\fi\else\ifE@\else\NESW@true\fi\fi
 \ifH@\else\ifV@\else
  \lams@slope
  \ifnum\lams@tani>\lams@tanii
   \lams@ht\ten@\p@\lams@wd\ten@\p@
   \multiply\lams@wd\lams@tanii\divide\lams@wd\lams@tani
  \else
   \lams@wd\ten@\p@\lams@ht\ten@\p@
   \divide\lams@ht\lams@tanii\multiply\lams@ht\lams@tani
  \fi
 \fi\fi
 \ifH@  \lams@harrow
 \else\ifV@ \lams@varrow
 \else \lams@darrow
 \fi\fi
}

\catcode`\@=\savecatcodeat
\let\savecatcodeat\undefined

\documentclass[tbtags,reqno]{amsart}        

\usepackage[tbtags]{amsmath}

\usepackage{amsthm,amsopn,amsfonts,amssymb,epsfig,pictex}            

\numberwithin{equation}{section}

\newtheorem{theorem}{Theorem}
\newtheorem{lemma}[theorem]{Lemma}
\newtheorem{proposition}[theorem]{Proposition}
\newtheorem{conjecture}[theorem]{Conjecture}

\newtheorem{definition}[theorem]{Definition}

\newtheorem{property}[theorem]{Property}
 
\theoremstyle{remark}

\newtheorem*{acknow}{\bf Acknowledgments}

\def\H{{\mathcal H}}

\def\endprf {\square}

\begin{document}

\title[Schur function analogs for a filtration of the symmetric function space]
{Schur function analogs for a filtration of the symmetric function space}

\author{L. Lapointe}
\address{ McGill University, 
Department of Mathematics and Statistics,  Montr\'eal, Qu\'ebec,
Canada, H3A 2K6}
\email{lapointe@math.mcgill.ca}

\author{J. Morse}
\address{University of Pennsylvania,
Department of Mathematics,
Philadelphia, PA 19104}
\email{morsej@math.upenn.edu}

\thanks{Research supported in part by NSF Grant \#0100179}

\begin{abstract}
We consider a filtration of the symmetric function space
given by $\Lambda^{(k)}_t$, the linear span of Hall-Littlewood 
polynomials indexed by partitions whose first part is not larger 
than $k$.  We introduce symmetric functions called the $k$-Schur 
functions, providing an analog for the Schur functions in the 
subspaces $\Lambda^{(k)}_t$.  We prove several properties for
the $k$-Schur functions including that they form a basis for 
these subspaces that reduces to the Schur basis when $k$ is large.  
We also show that the connection coefficients for the $k$-Schur 
function basis with the  Macdonald polynomials 
belonging to $\Lambda^{(k)}_t$ are polynomials in $q$ and $t$ with integral 
coefficients.  In fact, we conjecture that these integral
coefficients are actually positive, and give
several other conjectures generalizing Schur function theory.
\end{abstract}

\maketitle

\section{Introduction }

Let $\Lambda$ be the ring of symmetric functions
in the variables $x_1,x_2,\ldots$, with coefficients
in $\mathbb Q(q,t)$, for parameters $q$ and $t$. 
The Schur functions, $s_\lambda[X]$,  form a fundamental basis of 
$\Lambda$, with central roles in fields such as representation 
theory and algebraic geometry.  
For example, the Schur functions can be identified
with the characters of irreducible representations
of the symmetric group, and their products are equivalent to the
Pieri formulas for multiplying Schubert varieties in the intersection ring
of a Grassmannian.  
Furthermore, the connection coefficients of 
the Schur function basis with various bases
such as the homogeneous symmetric functions, 
the Hall-Littlewood polynomials, and the Macdonald polynomials,
are positive and have representation theoretic interpretations.
In the case of the Macdonald polynomials, $H_\lambda[X;q,t]$,
this expansion takes the form
\begin{equation} \label{MacKos}
H_{\lambda}[X;q,t] = \sum_{\mu} K_{\mu \lambda}(q,t) \, s_{\mu}[X] \, , \qquad
K_{\mu \lambda}(q,t) \in \mathbb N[q,t] \, ,
\end{equation}
where $K_{\mu \lambda}(q,t)$ are known as the $q,t$-Kostka polynomials. 
The representation theoretic interpretation for these polynomials is
given in \cite{[Ga],[Ga2],[Ha]}.

Here, we consider the filtration 
$\Lambda_t^{(1)}\subseteq \Lambda_t^{(2)}\subseteq
\cdots \subseteq\Lambda_t^{(\infty)}=\Lambda$,
given by the linear span of 
Hall-Littlewood polynomials indexed by $k$-bounded partitions.
That is,
\begin{equation}
\Lambda^{(k)}_t = 
{\mathcal L} \{H_{\lambda}[X;t] \}_{\lambda;\lambda_1 \leq k}\, , 
\qquad k=1,2,3,\dots \, .
\end{equation} 
We introduce a new family of symmetric functions
that are indexed by $k$-bounded partitions, denoted $s_\lambda^{(k)}[X;t]$, and prove 
a number of properties for these functions.
In particular, we show that they form bases for the 
subspaces, $\Lambda^{(k)}_t$.  Our functions will be called 
the $k$-Schur functions since they appear to play a role for  
$\Lambda^{(k)}_t$ that is analogous to the role of 
the Schur functions for $\Lambda$.  That is, the $k$-Schur functions 
give rise to the generalization of many Schur positivity properties.
Details are given following a brief outline of
our construction for $s_\lambda^{(k)}[X;t]$.

The characterization of $s_\lambda^{(k)}[X;t]$
relies on a $t$-generalization for Schur function products.
More precisely, for any partition sequence 
$S=(\lambda^{(1)},\dots,\lambda^{(\ell)})$,
a $t$-analog for the product 
$s_{\lambda^{(1)}}[X]\cdots s_{\lambda^{(\ell)}}[X]$
was introduced in \cite{[LLT],[Shi],[S1]}.
We first prove that a very particular subset of these
generalized products forms a basis for $\Lambda^{(k)}_t$.
The elements of our basis, denoted $G_\lambda^{(k)}[X;t]$,  
are the generalized Schur products with
sequence $S$ obtained by splitting $\lambda$ into pieces 
that depend on $k$.
$G_\lambda^{(k)}[X;t]$  
are thus called $k$-split polynomials.  
These polynomials
are essential in our definition for the $k$-Schur 
functions as we use a linear operator on $\Lambda^{(k)}_t$ 
defined by
\begin{equation}
{\bar T}_i^{(k)} \, G_{\lambda}^{(k)}[X;t] =
\begin{cases}
G_{\lambda}^{(k)}[X;t] & \text{ if } \lambda_1=i \\
0 & \text{ otherwise }
\end{cases} \, .
\end{equation}

The final ingredient needed to define the 
$k$-Schur functions is the vertex operator, 
$B_i$, introduced in \cite{[Ji]} to recursively
build the Hall-Littlewood polynomials.
More precisely,
\begin{equation} 
H_{\lambda_1,\ldots,\lambda_\ell}[X;t]
= B_{\lambda_1} \cdots B_{\lambda_\ell} \cdot 1 \, .
\end{equation}
Analogously to this relation,
we now define the $k$-Schur function 
for $k$-bounded $\lambda=(\lambda_1,\dots,\lambda_{\ell})$, by 
\begin{equation} \label{kkschur}
s_{\lambda}^{(k)}[X;t] = {\bar T}_{\lambda_1}^{(k)} B_{\lambda_1} \cdots
{\bar T}_{\lambda_\ell}^{(k)} B_{\lambda_\ell} \cdot 1 \, .
\end{equation}

Our work to characterize this basis was originally motivated
by two conjectures suggesting that the $k$-Schur functions 
play a central role in the understanding of the $q,t$-Kostka polynomials.  
Together, these conjectures refine relation  \eqref{MacKos}.
That is, for any $k$-bounded partition $\lambda$, 
\begin{eqnarray} \label{posischur}
i) &
s_{\lambda}^{(k)}[X;t]= \sum_{\mu \geq \lambda} v_{\mu \lambda}^{(k)}(t)\,
s_{\mu}[X] \, , \qquad v_{\mu \lambda}^{(k)}(t) \in \mathbb N[t] \, .
\label{con1}
\\
\nonumber \\ \label{posiMac}
ii) &
H_{\lambda}[X;q,t] = \sum_{\mu;\mu_1\leq k} 
K_{\mu \lambda}^{(k)}(q,t) \, s_{\mu}^{(k)}[X;t] \, , \qquad
K_{\mu \lambda}^{(k)}(q,t) \in \mathbb N[q,t] \, .
\label{con2}
\end{eqnarray}
Both conjectures hold when $k=2$ (Section~\ref{sec7})
and we shall prove that $v_{\mu\lambda}^{(k)}(t)\in\mathbb Z[t]$
and that $K_{\mu\lambda}^{(k)}(q,t)\in\mathbb Z[q,t]$
for all $k$.  Tables of coefficients $v_{\mu \lambda}^{(k)}(t)$ 
and $K_{\mu \lambda}^{(k)}(q,t)$ are included in section~\ref{tabs}
to illustrate these conjectures.  Our examples suggest 
a stronger property,
\begin{equation}
0\subseteq K_{\mu\lambda}^{(k)}(q,t)\subseteq
K_{\mu\lambda}(q,t)
\end{equation}
where for two polynomials $P,Q\in\mathbb Z[t]$, 
$P\subseteq Q$ means $Q-P\in \mathbb N[q,t]$.

More generally, it develops that properties of the $k$-Schur functions, 
with a number of conjectures, provide a $k$-generalization for
the properties that make the Schur functions important to the 
theory of symmetric functions.
In particular, we prove that the $s_\lambda^{(k)}[X;t]$ form a basis for 
$\Lambda^{(k)}_t$ and that the $k$-Schur functions of 
$\Lambda^{(\infty)}_t=\Lambda$ are indeed the Schur functions 
themselves. 
Conjectural evidence for the significance of the $k$-Schur functions
includes a $k$-analog of partition conjugation
and generalizations of Pieri and Littlewood-Richardson rules. 
Consequently, in the case of the multiplicative action of $h_1[X]$, 
a $k$-analog of the Young Lattice is induced.
Further, we have observed that the $k$-Schur functions, expanded
in terms of $k$-Schur functions in two sets of variables, 
have coefficients in $\mathbb N[t]$. This is 
a special property of Schur functions that is not shared 
by the Hall-Littlewood or Macdonald functions.  
Finally, the $k$-Schur functions
of $\Lambda^{(k)}_t$, when embedded in  
$\Lambda^{(k')}_t$ 
for $k'>k$, seem to decompose positively in terms of $k'$-Schur functions:
\begin{equation}
s_{\lambda}^{(k)}[X;t]=s_{\lambda}^{(k')}[X;t]+ \sum_{\mu > \lambda}
v_{\mu \lambda}^{(k \to k')}(t) \, s_{\mu}^{(k')}[X;t] \, ,\quad {\rm{where~}}
v_{\mu \lambda}^{(k \to k')}(t) \in \mathbb N[t] \, .
\end{equation}

It happens, remarkably, that not all of the $k$-Schur functions 
need to be constructed using \eqref{kkschur}.  
For each $k$, there is a subset of $s_\lambda^{(k)}[X;t]$, 
called the irreducible $k$-Schur functions, from which all other 
$k$-Schur functions may be constructed 
\cite{[LM2]} by 
simply applying a succession of certain operators.
The elements of this set are the $k$-Schur functions 
indexed by partitions with no more than $i$ parts equal to $k-i$,
and the operators are vertex operators \cite{[ZS]} associated to rectangularly 
shaped partitions $(\ell^{k+1-\ell})$ for $\ell=1,\dots,k$.
That is,
\begin{equation}
s^{(k)}_\lambda[X;t]
=
t^c \, 
 B_{R_1}
 B_{R_2}
\ldots
 B_{R_\ell} \, 
 s^{(k)}_\mu[X;t] \quad\text{with}\quad
c \in\mathbb N \, ,
\label{atirr}
\end{equation}
for an irreducible $s_\mu^{(k)}$,
and vertex operators $B_R$, where
$R$ is a partition of rectangular shape.

\medskip

Since the Hall-Littlewood polynomials at $t=1$ are the 
complete symmetric functions
\begin{equation}
H_{\lambda}[X;1] = 
h_{\lambda_1}[X]\, h_{\lambda_2}[X] \cdots h_{\lambda_\ell}[X]
\,,
\end{equation}
we see that $\Lambda^{(k)}_t$ reduces to the polynomial
ring $\Lambda^{(k)}=\mathbb Q[h_1,\ldots,h_k]$.
Each of the properties held by the $k$-Schur functions 
has a specialization in this subring.
In particular, since $B_R$ is simply multiplication 
by the Schur function $s_R$ 
when $t=1$, relation \eqref{atirr} reduces to
\begin{equation}
s^{(k)}_\lambda[X]
=
s_{R_1}[X]
s_{R_2}[X]
\ldots
s_{R_\ell}[X]
s^{(k)}_\mu[X]
\,.
\label{surmul}
\end{equation}
The irreducible $k$-Schur functions 
thus constitute a natural basis
for the quotient ring  $\Lambda^{(k)}/{\mathcal I}_k$, where
${\mathcal I}_k$ is the ideal generated by 
Schur functions indexed by partitions of the form $(\ell^{k+1-\ell})$.

This work is closely related to \cite{[LLM]}, 
where symmetric functions called atoms are 
introduced in terms of tableaux. 
We conjecture that the $k$-Schur function 
$s_\lambda^{(k)}[X;t]$ defined here, 
is equivalent to the atom $A_\lambda^{(k)}[X;t]$ of 
\cite{[LLM]}\footnote
{The change in terminology is to prevent confusion 
with other objects know as atoms in tableaux combinatorics.}.
Thus, a combinatorial interpretation for the $k$-Schur 
functions can be found in \cite{[LLM]}, as well as many 
other properties that our functions appear to obey.  
In particular, a discussion of a $k$-generalization for the Young 
lattice is present.  The reader should be aware of notational differences 
in \cite{[LLM]} such as 
$A_{\lambda}^{(k)}[X;t] \leftrightarrow s_{\lambda}^{(k)}[X;t]$,
$V_k \leftrightarrow \Lambda^{(k)}_t$, and $S_{\lambda}[X] \leftrightarrow
s_{\lambda}[X]$.

The article proceeds as follows:
after the introduction of basic definitions in Section~\ref{sec1}, 
the case $t=1$ is addressed in 
Sections~\ref{sec3}, \ref{secirred} and \ref{sec5}.
Although some of what appears in these sections can be obtained 
from the results in Sections~\ref{sec6} and \ref{sec7}, 
owing to its greater simplicity, we find it pedagogical 
to treat the case $t=1$ independently.  
In Section~\ref{sec6}, an important
unitriangularity property of
the generalized Schur products is used.  However, 
the proof has been relegated to the Appendix, where 
known properties of these products including a
Morris-type recursion can also be found.

\smallskip

\begin{acknow}
{\it
The enthusiasm from A. Garsia and A. Lascoux greatly contributed
to this work and we are thankful to M. Zabrocki
for helping us with \cite{[S1]}.  L. Lapointe thanks
L. Vinet for his support.  J. Morse held an NSF grant
for part of the period devoted to this research.
ACE \cite{[V]} was instrumental towards this work.
}
\end{acknow}

\smallskip

\section{Definitions} \label{sec1}

\subsection{Partitions}
Symmetric polynomials are indexed by partitions, sequences of
non-negative integers $\lambda =(\lambda_1,\lambda_2,\ldots)$ with
$\lambda_1 \ge \lambda_2 \ge \dots$.  
The number of non-zero parts in $\lambda$ is denoted $\ell(\lambda)$
and the degree of $\lambda$ is $|\lambda| = \lambda_1 + \dots + 
\lambda_{\ell(\lambda)}$.  We say that $\lambda$ is a partition of $n$,
denoted $\lambda \vdash n$, if $|\lambda|=n$.
We use the dominance order on partitions with $|\lambda|=|\mu|$, where
$\lambda\leq\mu$ when $\lambda_1+\cdots+\lambda_i\leq
\mu_1+\cdots+\mu_i$ for all $i$.  
Given two partitions $\lambda$ and $\mu$,
$\lambda \cup \mu$ stands for the partition 
rearrangement of the parts of $\lambda$ and $\mu$. 
Note that if $\lambda \leq \mu$ and $\nu \leq \omega$, then
$\lambda \cup \nu \leq \mu \cup \omega$.  

Any partition $\lambda$ has an associated Ferrers diagram 
with $\lambda_i$ lattice squares in the $i^{th}$ row, 
from the bottom to top.  For example,
\begin{equation}
\lambda\,=\,(4,2)\,=\,
{\tiny{\tableau*[scY]{ & \cr & & & }}} \, .
\end{equation}
For each cell $s=(i,j)$ in the diagram of $\lambda$, let
$\ell'(s), \ell(s), a(s)$ and $a'(s)$ be respectively the number of
cells in the diagram of $\lambda$ to the south, north, east and west
of the cell $s$.  The hook-length of any cell in $\lambda$, is
$h_s(\lambda)=\ell(s)+a(s)+1$.  In the example,
$h_{(1,2)}(4,2)=2+1+1$.  The {\it main hook-length} of $\lambda$, $h_M(\lambda)$,
is the hook-length of the cell $s=(1,1)$ in the diagram of $\lambda$.  
Therefore, $h_{M}\bigl((4,2)\bigr)=5$.

The conjugate $\lambda'$ of a  partition $\lambda$ is defined
by the reflection of the Ferrers diagram about the main diagonal.
For example, the conjugate of (4,2) is
\begin{equation}
\lambda' \,=\,
{\tiny{\tableau*[scY]{ \cr \cr & \cr &   }}}
\,=\,(2,2,1,1)\,.
\end{equation}
A skew shape $\lambda/\mu$ is the diagram obtained by 
deleting the diagram of $\mu$ from $\lambda$.  
\begin{equation}
(6,5,4,3,1)/(4,2) = 
{\tiny{\tableau*[scY]
{ \cr & & \cr & & & 
\cr \bl &\bl  & &  &
\cr \bl &\bl  & \bl &\bl  & & 
}}}
\end{equation}

A partition $\lambda$ is said to be $k$-{\it bounded} if its first part is 
not larger than $k$, i.e,  if $\lambda_1 \leq k$. 
We associate to any $k$-bounded partition $\lambda$
a sequence of partitions, 
$\lambda^{\to k}= (\lambda^{(1)},\lambda^{(2)},\ldots,\lambda^{(r)})$,
called the $k$-{\it split} of $\lambda$.
$\lambda^{\to k}$ is obtained by partitioning $\lambda$ 
(without rearranging the entries) into partitions $\lambda^{(i)}$
where $h_M(\lambda^{(i)})= k$, for all $i<r $.  For example,
$(3,2,2,2,1,1)^{\to 3}=\bigl((3),(2,2),(2,1),(1)\bigr)$ 
and $(3,2,2,2,1,1)^{\to 4}=\bigl((3,2),(2,2,1),(1)\bigr)$ .
Equivalently, the diagram of $\lambda$ is cut
horizontally into partitions with main hook-length $k$.
\begin{equation}
{\tiny{\tableau*[scY]{ \cr \cr & \cr & \cr & \cr & & \cr}}} \quad
\begin{matrix}
& {\tiny{\tableau*[sbY]{  \cr}}} \hfill \\
& {\tiny{\tableau*[sbY]{  \cr & \cr}}} \hfill \\
\longrightarrow^{(3)} & {\tiny{\tableau*[sbY]{ & \cr  & \cr}}} \hfill \\
& {\tiny{\tableau*[sbY]{  & & \cr}}} \hfill
\end{matrix}
\qquad
\text{and}
\qquad
{\tiny{\tableau*[scY]{ \cr \cr & \cr & \cr & \cr & & \cr}}} \quad
\begin{matrix}
& {\tiny{\tableau*[sbY]{  \cr}}} \hfill \\
\longrightarrow^{(4)} & {\tiny{\tableau*[sbY]{\cr  & \cr  & \cr}}} \hfill \\
& {\tiny{\tableau*[sbY]{ & \cr & & \cr}}} \hfill
\end{matrix}
\, .
\end{equation}
The last partition in the sequence $\lambda^{\to k}$ may have main hook-
length less than $k$.
It is important to note that $\lambda^{\to k}=(\lambda)$
when $h_M(\lambda)\leq k$.

\subsection{Symmetric functions} The power sum $p_i(x_1,x_2,\ldots)$ is
\begin{equation}
p_i(x_1,x_2,\ldots) = x_1^i+x_2^i+\cdots \, ,
\end{equation}
and for a partition $\lambda=(\lambda_1,\lambda_2,\dots)$, 
\begin{equation}
p_{\lambda}(x_1,x_2,\ldots)=
p_{\lambda_1}(x_1,x_2,\ldots)\, p_{\lambda_2}(x_1,x_2,\ldots)\cdots \,.
\end{equation}
We employ the notation of $\lambda$-rings, needing only the formal 
ring of symmetric functions $\Lambda$ to act on the ring of rational 
functions in $x_1,\dots,x_N,q,t$, with coefficients in $\mathbb R$.
The action of a power sum $p_i$ on a rational function is, by definition,
\begin{equation}
p_{i} \left[ \frac{\sum_{\alpha} c_{\alpha} u_{\alpha} }
{ \sum_{\beta} d_{\beta} v_{\beta} } \right]
 =\frac{\sum_{\alpha} c_{\alpha} u_{\alpha}^i }
{ \sum_{\beta} d_{\beta} v_{\beta}^i},
\label{actp}
\end{equation}
with $c_{\alpha},d_{\beta} \in \mathbb R$ and $u_{\alpha}, v_{\beta}$
monomials in $x_1,\dots,x_N,q,t$.  Since the power sums form a basis
of the ring $\Lambda$, any symmetric function has a 
unique expression in terms of power sums, and \eqref{actp} extends to an 
action of $\Lambda$ on rational functions.
In particular $f[X]$, the action of a symmetric function
$f$ on the monomial $X=x_1+\cdots+x_N$, is simply $f(x_1,\ldots,x_N)$.
In the remainder of the article, we will always consider the number of
variables $N$ to be infinite, unless otherwise specified.

\smallskip

The complete symmetric function $h_r[X]$ is
\begin{equation}
h_r[X] = \sum_{1\leq i_1\leq i_2\leq\cdots\leq i_r}
x_{i_1} x_{i_2}\cdots x_{i_r} \, ,
\end{equation}
and $h_{\lambda}[X]$ stands for the homogeneous symmetric function
\begin{equation}
h_{\lambda}[X]=
h_{\lambda_1}[X]\, h_{\lambda_2}[X]\cdots \,.
\end{equation} 
In the same way, the elementary symmetric function $e_r[X]$ is
\begin{equation}
e_r[X] = \sum_{1\leq i_1< i_2<\cdots< i_r}
x_{i_1} x_{i_2}\cdots x_{i_r} \, ,
\end{equation}
and $e_{\lambda}[X]$ stands for
\begin{equation}
e_{\lambda}[X]=
e_{\lambda_1}[X]\, e_{\lambda_2}[X]\cdots \,.
\end{equation} 

Although the Schur functions may be characterized in many ways, 
here it will be convenient to use the Jacobi-Trudi 
determinantal expression:
\begin{equation}
s_\lambda[X] = \text{det}\Bigl | h_{\lambda_i+j-1}[X] 
\Bigr|_{1 \leq i,j \leq \ell(\lambda)}=
\text{det}\Bigl | e_{\lambda_i'+j-1}[X] \Bigr|_{1 \leq i,j \leq \ell(\lambda)} \, ,
\label{jt}
\end{equation}
where $h_r[X]=e_r[X]=0$ if $r < 0$.  
Note, in particular, $s_r[X]=h_r[X]$ and $s_{1^r}[X]=e_r[X]$. 

The homomorphism $\omega$, which is an involution on $\Lambda$, 
is defined by
\begin{equation} \label{homomor}
\omega(h_r[X])= e_r[X] \, ,
\end{equation}
and is such that $\omega(s_{\lambda}[X])=s_{\lambda'}[X]$.

We recall that the Macdonald scalar product,  
$\langle \ , \ \rangle_{q,t}$, 
on $\Lambda \otimes \mathbb Q(q,t)$ is defined by setting 
\begin{equation}
\langle p_\lambda[X], p_\mu[X] \rangle_{q,t}
        =\delta_{\lambda \mu } \, z_\lambda \prod_{i=1}^{\ell(\lambda)} \frac{1-q^{\lambda_i}}{1-t^{\lambda_i}} \, ,
\end{equation}
where for a partition $\lambda$ with $m_i(\lambda)$ parts equal to $i$,
we associate the number
\begin{equation} \label{1}
z_\lambda
        = 1^{m_1} m_1 !  \, 2^{m_2} m_2! \dotsm
\end{equation}
When $q=t$, this scalar product does not depend
on a parameter and then satisfies
\begin{equation} \label{schurscalar}
\langle s_{\lambda}[X],s_{\mu}[X]\rangle= \delta_{\lambda \mu} \, .
\end{equation}

The Macdonald integral forms 
$J_\lambda [X; q,t]$ are uniquely characterized \cite{[Ma]} by
\begin{align}
  \mathrm{(i)} \ &  \langle J_\lambda, J_\mu \rangle_{q,t} = 0, \qquad \text{if } \lambda \ne \mu , \\
  \mathrm{(ii)} \ &  J_\lambda[X;q,t] = \sum_{\mu \le \lambda} v_{\lambda\mu }(q,t) 
s_\mu[X] \, , 
\quad\text{with}\quad v_{\lambda \mu}(q,t) \in \mathbb Q(q,t)\,,  \\
  \mathrm{(iii)} \ & v_{\lambda\lambda}(q,t)= 
\prod_{s \in \lambda} (1-q^{a(s)} t^{\ell (s)+1}) \, .
\end{align}
Here, we use a modification of the Macdonald 
integral forms that is obtained by setting
\begin{equation}
H_{\lambda}[X;q,t]=J_{\lambda}[X/(1-t);q,t]= \sum_{\mu} K_{\mu \lambda}(q,t)\, 
s_{\mu}[X] \, ,
\end{equation}
with the coefficients $K_{\mu \lambda}(q,t) \in \mathbb N[q,t]$ 
known as the $q,t$-Kostka polynomials.
In the case $q=0$, $J_{\lambda}[X;q,t]$ reduces to the Hall-Littlewood
polynomial, $J_{\lambda}[X;q,t]=Q_{\lambda}[X;t]$.  Again, we shall
use a modification of the Hall-Littlewood polynomials,
\begin{equation}
H_{\lambda}[X;t]=H_{\lambda}[X;0,t]=Q_{\lambda}[X/(1-t);t]= 
s_{\lambda}[X]+\sum_{\mu> \lambda} K_{\mu \lambda}(t)\, 
s_{\mu}[X] \, ,
\label{HallinS}
\end{equation}
with the coefficients $K_{\mu \lambda}(t) \in \mathbb N[t]$ known 
as the Kostka-Foulkes polynomials.
Then, in the limit $t=1$, we have  $H_{\lambda}[X;1]=h_{\lambda}[X]$, giving
\begin{equation} \label{kostka}
h_{\lambda}[X]=H_{\lambda}[X;1]= 
s_{\lambda}[X]+\sum_{\mu>\lambda} K_{\mu \lambda}\, 
s_{\mu}[X] \, ,
\end{equation}
with the coefficients $K_{\mu \lambda} \in \mathbb N$ known 
as the Kostka numbers.

For partitions $\lambda$ and $\mu$, we have
\begin{equation}
\label{LRprod}
s_{\lambda}[X] \, s_{\mu}[X] = 
\sum_{{\rho : |\rho|=|\lambda|+|\mu|}}
c_{\lambda \mu}^{\rho} \, s_{\rho}[X] \, ,
\end{equation}
where the coefficients $ c_{\lambda \mu}^{\rho} \in \mathbb N$
are the Littlewood-Richardson coefficients.  They satisfy
\begin{property}  
$c_{\lambda \mu}^{\lambda \cup \mu}=1$, and
$c_{\lambda \mu}^{\rho}=0$ unless $\lambda \cup \mu \leq \rho$ and
$\left(\lambda' \cup \mu' \right)'\geq \rho$. 
\label{LR}
\end{property}
\noindent 
{\bf Proof.}\quad The unitriangular relation (\ref{kostka}) gives
\begin{equation}
s_{\lambda}[X]= h_{\lambda}[X]+\sum_{\mu>\lambda} K_{\mu \lambda}^{-1}\, 
h_{\mu}[X] \, ,
\end{equation}
which implies
\begin{equation}
s_{\lambda}[X]\, s_{\mu}[X]= \sum_{
\nu \geq \lambda; \gamma \geq
\mu } K_{\nu \lambda}^{-1} K_{\gamma \mu}^{-1} \, h_{\nu \cup \gamma}[X]
= h_{\lambda \cup \mu} + \sum_{\rho > \lambda \cup \mu} d^{~\rho}_{\lambda \mu}\, 
h_{\rho} \, ,
\end{equation}
since $\nu \geq \lambda $ and $\gamma \geq \mu$ imply $\nu \cup \gamma \geq 
\lambda \cup \mu$.  Then, from (\ref{kostka}),  
$c_{\lambda \mu}^{\lambda \cup \mu}=1$ and
$c_{\lambda \mu}^{\rho}=0$ unless $\lambda \cup \mu \leq \rho$.
Similarly, it can be shown that
$\left(\lambda' \cup \mu' \right)'\geq \rho$
by applying the involution $\omega$ to (\ref{kostka}).

The Pieri rules are the cases
$\mu=(r)$ and $\mu=(1^r)$ of expression (\ref{LRprod}).
Then
\begin{equation}
h_r[X] s_{\lambda}[X]  = \sum_{\mu} s_{\mu}[X] \quad
\quad
\text{and}
\quad
\quad
e_r[X] s_{\lambda}[X]  = \sum_{\nu} s_{\nu}[X] \, , 
\end{equation}
where the sums run over all $\mu$'s sand $\nu$'s such that $\mu/\lambda$ and
$\nu/\lambda$ are respectively a horizontal $r$-strip and a vertical
$r$-strip.

\section{The $k$-Schur functions} \label{sec3}
Here, and in the following section, we study the space
$\Lambda_{t}^{(k)}$ when $t=1$.  In this case,
the subspace reduces to a subring of $\Lambda$, defined
by
\begin{equation} \label{ringhk}
\Lambda^{(k)}=\Lambda_{t=1}^{(k)}=  
{\mathcal L} \left\{ h_\lambda[X] \right\}_{\lambda;\lambda_1\leq k}=
{\mathcal L} \left\{ e_\lambda[X] \right\}_{\lambda;\lambda_1\leq k} \, .
\end{equation}
The last equality holds since the determinantal expression 
\eqref{jt} of $s_r[X]=h_r[X]$, in terms of
elementary functions, gives only terms of the type $e_i[X]$, 
where $i\leq r$.

\subsection{$k$-split polynomials}
Our construction relies on the introduction of a 
family of symmetric polynomials called the $k$-split 
polynomials.  If $\lambda^{\to k} = 
(\lambda^{(1)},\lambda^{(2)},\ldots,\lambda^{(n)})$
is the $k$-split of a $k$-bounded partition $\lambda$,
then we define the $k$-split polynomial by
\begin{equation}
G_{\lambda}^{(k)}[X] = 
s_{\lambda^{(1)}}[X]
\,
s_{\lambda^{(2)}}[X]
\cdots
s_{\lambda^{(n)}}[X] \, .
\end{equation}
A number of properties can be derived from this definition.
For example, because $\lambda^{\to k}=(\lambda)$ 
when $h_M(\lambda)\leq k$, it follows immediately that

\medskip

\begin{property}
\begin{equation}
G_\lambda^{(k)}[X] = s_\lambda[X]\, ,\qquad\text{for}\quad h_M(\lambda)\leq k
\,.
\end{equation}
\label{kpol1}
\end{property}
\noindent Moreover, since $G_\lambda^{(k)}[X]$ is the product 
$s_{\lambda^{(1)}}[X] \cdots s_{\lambda^{(n)}}[X]$,
Property~\ref{LR} of the Littlewood-Richardson rule 
applied repeatedly implies that 
$\lambda= \lambda^{(1)}\cup \cdots \cup \lambda^{(n)}$ must index
the minimal element in the Schur function expansion 
of $G_\lambda^{(k)}[X]$.  That is,

\medskip

\begin{property}
For any $k$-bounded partition $\lambda$,
\begin{equation}
G_\lambda^{(k)}[X]  = s_\lambda[X] + \sum_{\mu>\lambda} 
a_{\lambda \mu}^{(k)}\,s_\mu[X] \, ,\qquad 
a_{\lambda \mu}^{(k)} \in \mathbb N \, .
\end{equation}
\label{kpol2}
\end{property}
\noindent
Finally, 
since $(k,\lambda)^{\to k} =\bigl((k),\lambda^{\to k}\bigr)$
where $(k,\lambda)=(k,\lambda_1,\lambda_2,\dots)$,
we have

\medskip

\begin{property}
For any $k$-bounded partition $\lambda$,
\begin{equation}
s_k[X]\,G_\lambda^{(k)}[X] = G_{(k,\lambda)}^{(k)}[X] 
\, .
\end{equation}
\label{kpol3}
\end{property}

In fact, it happens that these polynomials form a basis for
the subring $\Lambda^{(k)}$.

\medskip

\begin{theorem} \label{theobaseG}
$\{ G_\lambda^{(k)}[X] \}_{\lambda_1\leq k}$ forms a basis for 
$\Lambda^{(k)}$.
\end{theorem}
\noindent 
{\bf Proof.}\quad  
$\Lambda^{(k)}$ is the linear span of homogeneous symmetric functions
over all $k$-bounded partitions.  
Since the elements $G_\lambda^{(k)}[X]$ are also indexed by $k$-bounded
partitions and by Property~\ref{kpol2},
are linearly independent, 
it suffices to show that these elements lie in $\Lambda^{(k)}$.  
That is, for a $k$-bounded partition $\lambda$, 
we must be able to expand $G_\lambda^{(k)}[X]$ in terms
of $h_{\mu_1}[X]h_{\mu_2}[X]\cdots$ where $\mu_i\leq k$ for all $i$.
Observe that any entry $h_j[X]$ in the
determinantal expression \eqref{jt} for $s_\nu$
can be indexed by $j$ no larger than
$\nu_1+\ell(\nu)-1=h_M(\nu)$.  
Since $G_\lambda^{(k)}[X]=s_{\lambda^{(1)}}[X]s_{\lambda^{(2)}}[X]\cdots$
with $h_M(\lambda^{(i)})\leq k$ for all $i$,
$G_\lambda^{(k)}[X]$ is a product of 
determinants having entries $h_j[X]$, with $j\leq k$.
\hfill $\square$

\medskip

We now have that $G_\lambda^{(k)}\in \Lambda^{(k)}$
for all $k$-bounded partitions $\lambda$.
Therefore, these polynomials can be expanded in terms of
$h_\mu$ with $\mu_1\leq k$.  In fact, 
since $G_\lambda^{(k)}$ is unitriangularly 
related to $s_\lambda$ by Property~\ref{kpol2},
and $s_\lambda$ is unitriangularly related
to $h_\lambda$ \eqref{kostka}, we have 

\medskip

\begin{property}
\begin{equation}
G_{\lambda}^{(k)}[X] = h_\lambda[X] + \sum_{\mu>\lambda\atop\mu_1\leq k} 
g_{\lambda \mu}^{(k)}\,h_\mu[X] \, ,\qquad 
g_{\lambda \mu}^{(k)} \in \mathbb Z \, .
\label{uniGh}
\end{equation}
\end{property}

\medskip

\subsection{$k$-Schur functions}
Although the $k$-split polynomials form a basis for $\Lambda^{(k)}$, 
they do not have the fundamental role for $\Lambda^{(k)}$
that the Schur functions have for $\Lambda$.  
For example, since
$$
s_2[X]\,  G_{3,1,1}^{(3)}[X] \,=\,
G_{3,2,1,1}^{(3)}[X]- G_{3,2,2}^{(3)}[X]
$$
reveals a negative coefficient, the $k$-split basis 
does not satisfy a refinement of the Pieri Rule.  
Additionally, there are no partitions $\mu$ such that
the involution $\omega$ sends
$G_{\lambda}^{(k)}[X]$ to a single element $G_\mu^{(k)}[X]$
as with the Schur functions. 
However, these polynomials do play a key role in the 
construction of our Schur analog, $s_\lambda^{(k)}[X]$.  
It happens that we use a projection operator, $T_j^{(k)}$, 
acting linearly on $\Lambda^{(k)}_t$ and defined by
\begin{equation}
T_j^{(k)} \, G_\lambda^{(k)}[X] = 
\begin{cases}
G_\lambda^{(k)}[X] & \text{if } \lambda_1 = j\\
0 & \text{otherwise }
\end{cases}
\, .
\end{equation}
Note that this projection operator
sends an element of $\Lambda^{(k)}$ into an element belonging to
the linear span of $k$-split polynomials
whose first parts all equal to $j$. 
We can now introduce the $k$-Schur functions,
our fundamental objects in $\Lambda^{(k)}$.

\begin{definition} \label{kschurdef}
For a $k$-bounded partition $\lambda$, 
the $k$-Schur function is defined recursively by
\begin{equation}
s_\lambda^{(k)}[X] = T_{\lambda_1}^{(k)}\,s_{\lambda_1}[X]\,
s_{\lambda_2,\lambda_3,\dots}^{(k)}[X] \, , 
\quad
\text{with} \quad s_{()}^{(k)}[X]=1\,.
\end{equation}
\end{definition}
Explicit examples of the expansion of $k$-Schur 
functions in terms of Schur functions are given,
letting $t=1$, in the tables of Subsection~\ref{tabkschur}.

\medskip

Since the space $\Lambda^{(k)}$ is simply $\Lambda$
when $k\to\infty$, we expect to recover 
$s_\lambda$ from $s_\lambda^{(k)}$ in this case. 
The following property supports our assertion that 
the $k$-Schur functions generalize the Schur functions.

\medskip

\begin{property} \label{largek}
For any $k$-bounded partition $\lambda$,
\begin{equation}
s_\lambda^{(k)}[X] = s_\lambda[X] \, ,\qquad\text{for}\quad 
h_M(\lambda)\leq k \,.
\end{equation}
\end{property}
\noindent{\bf Proof.}\quad
Let $\lambda=(\lambda_1,\lambda_2\ldots,)$
be a partition with $\lambda_1\leq k$
and $h_M(\lambda)\leq k$.  
In the case that $\lambda=()$,
by definition, $s_{()}^{(k)}[X]=1=s_{()}[X]$.
Let $\hat \lambda=(\lambda_2,\lambda_3,\dots)$
and assume by induction on $\ell(\lambda)$, that 
$s_{\hat \lambda}^{(k)}[X] =s_{\hat \lambda}[X]$.
We proceed to show $s_\lambda^{(k)}[X]=s_\lambda[X]$.
From \eqref{jt}, we have that
\begin{equation}
s_{\lambda}[X] = 
\det \left|\begin{matrix}
h_{\lambda_1}[X] & h_{\lambda_1+1}[X] & \cdots &  h_{\lambda_1+\ell-1}[X] \\
h_{\lambda_2-1}[X] & h_{\lambda_2}[X] & \cdots & h_{\lambda_2+\ell-2}[X] \\
\vdots &  \vdots & \ddots & \vdots \\
h_{\lambda_\ell-\ell+1}[X] & 
h_{\lambda_{\ell}-\ell+2}[X] & \cdots & h_{\lambda_\ell}[X] 
 \end{matrix} \right | \, . 
\end{equation}
The expansion of this determinant about the first row
gives that the coefficient of $h_{\lambda_1}[X]$ 
is exactly $s_{\hat\lambda}[X]$.  
Further, since all other terms in the first row 
are of the form $h_{i}[X]$, with $i>\lambda_1$, we have
\begin{equation}
s_{\lambda}[X] = h_{\lambda_1}[X]s_{\hat \lambda}[X] + 
\sum_{\mu;\mu_1>\lambda_1} c_{\mu} \, h_{\mu}[X] \, .
\label{dum1}
\end{equation}
Property~\ref{kpol1} gives that $s_{\lambda}[X]=
G^{(k)}_{\lambda}[X]$ since $h_M(\lambda)\leq k$. 
Noting also that $s_{\lambda_1}=h_{\lambda_1}$, 
\eqref{dum1} becomes 
\begin{equation}
s_{\lambda_1}s_{\hat\lambda}[X] 
=
G^{(k)}_{\lambda}[X]  
- \sum_{\mu;\mu_1>\lambda_1} c_{\mu} \, h_{\mu}[X] 
\, .
\end{equation}
Since $s_{\lambda_1}[X]s_{\hat \lambda}[X]$
and $G_\lambda^{(k)}$ both belong to $\Lambda^{(k)}$,
$\sum_{\mu;\mu_1>\lambda_1} c_{\mu}h_{\mu}[X]$ is in
$\Lambda^{(k)}$ and can therefore be expanded in terms 
of $G_\mu^{(k)}$.  This expansion, by the 
unitriangularity property \eqref{uniGh}, yields
\begin{equation}
s_{\lambda_1}[X]s_{\hat \lambda}[X] = G_{\lambda}^{(k)}[X]- 
\sum_{\mu;\mu_1>\lambda_1} d_{\mu} \, G_{\mu}^{(k)}[X] \, .
\label{dum2}
\end{equation}
Since $s_\lambda^{(k)}[X] = T_{\lambda_1}^{(k)}
s_{\lambda_1}[X]s_{\hat \lambda}[X]$
by the induction hypothesis, applying $T_{\lambda_1}$ to 
this expression gives
\begin{equation}
s_\lambda^{(k)}[X] =
T_{\lambda_1}^{(k)} \left(G_{\lambda}^{(k)}[X]- 
\sum_{\mu;\mu_1>\lambda_1} d_{\mu} \, G_{\mu}^{(k)}[X]  \right) 
= G_{\lambda}^{(k)}[X]=s_{\lambda}[X] \, ,
\end{equation}
thus proving the property.
\hfill $\square$

It also happens that there is a unitriangularity property
satisfied by the $k$-Schur functions.

\medskip

\begin{property} \label{unish1}
For any $k$-bounded partition $\lambda$,
\begin{equation}
s_\lambda^{(k)}[X] = 
h_\lambda[X] + \sum_{\mu>\lambda\atop\mu_1 \leq k} 
d_{\lambda\mu}^{(k)}\,h_\mu[X]
\, , \qquad {\rm with } \quad d_{\lambda \mu}^{(k)} \in \mathbb Z \, .
\end{equation}
\label{tria}
\end{property}
\noindent{\bf Proof.}\quad
Let $\lambda$ be a $k$-bounded partition.
If $\lambda=(r)$ then
$h_M(\lambda)\leq k$ and
Property~\ref{largek} gives 
\begin{equation}
s_{r}^{(k)}[X]= s_{r}[X] = h_r[X]
\,.
\end{equation}
Let $\hat\lambda$ denote the partition $\lambda$ without its
first part and 
assume by induction on $\ell(\lambda)$ that
\begin{equation}
s_{\hat \lambda}^{(k)}[X]
=h_{\hat \lambda}[X] +
\sum_{\gamma>\hat \lambda\atop \gamma_1 \leq k} 
d_\gamma \, h_{\gamma}[X]
\, , \qquad {\rm with } \quad d_{\gamma} \in \mathbb Z \, .
\end{equation}
Substituting this expression into the definition
$s_\lambda^{(k)}[X] = T_{\lambda_1}s_{\lambda_1}
s_{\hat\lambda}^{(k)}$ gives
\begin{equation}
\begin{split}
s_\lambda^{(k)}[X] & = 
T_{\lambda_1}^{(k)} \,
\left(
h_{\lambda}[X]+
\sum_{\gamma>\hat \lambda; \gamma_1 \leq k} 
d_\gamma \, h_{(\lambda_1)\cup \gamma }[X] \right) \\
& = T_{\lambda_1}^{(k)} \,
\left(
h_{\lambda}[X] +
\sum_{\nu>\lambda; \nu_1 \leq k} 
c_\nu \, h_{\nu}[X] \right) 
\, , \qquad {\rm with } \quad c_{\nu} \in \mathbb Z \, .
\end{split}
\end{equation}
since $\gamma > \hat \lambda
\implies
(\lambda_1) \cup \gamma > (\lambda_1) \cup \hat \lambda=\lambda$.
The unitriangularity relation \eqref{uniGh} then gives
\begin{eqnarray} 
s_\lambda^{(k)}[X] & = &
T_{\lambda_1} \, 
\left(
G_\lambda^{(k)}[X] + \sum_{\mu>\lambda} v_\mu \, G_\mu^{(k)}[X]
\right)
\\ 
& = &
G_\lambda^{(k)}[X] + 
\sum_{\mu>\lambda ; \mu_1=\lambda_1}
v_\mu \, G_\mu^{(k)}[X]
\, , \qquad {\rm with } \quad v_{\mu} \in \mathbb Z \, .
 \label{unisG1proof}
\end{eqnarray}
Using formula \eqref{uniGh} again 
then proves the claim. 
\hfill $\square$

We now have the tools to prove that in fact,
the $k$-Schur functions do form a basis.

\medskip

\begin{theorem} 
The $k$-Schur functions form a basis for $\Lambda^{(k)}$.
That is
\begin{equation}
\Lambda^{(k)}
\, = \,
{\mathcal L} \left\{
s_\lambda^{(k)}[X]
\right\}_{\lambda_1\leq k}
\, .
\end{equation}
\end{theorem}
\noindent{\bf Proof.}
Since $\Lambda^{(k)}=
{\mathcal L} \{h_{\lambda}[X] \}_{\lambda;\lambda_1 \leq k}\, $,
the theorem follows immediately from the unitriangularity relation found
in Property~\ref{unish1}.
\hfill$\square$

\medskip

Property \ref{kpol3} of the $k$-split polynomials states that
$s_k[X] \,  G_\lambda^{(k)}[X]=G_{k,\lambda}^{(k)}[X]$.  
In the next section, we prove a similar property for the 
$k$-Schur functions and generalize the phenomenon.

\section{Irreducibility} \label{secirred}

Here we study the product of certain Schur functions
with a $k$-Schur function.  In particular, we 
are concerned with partitions of the form $(\ell^{k+1-\ell})$
for $\ell=1,\ldots,k$.  Hereafter, such a partition is 
referred to as a $k$-rectangle.  A $k$-bounded partition with
no more than $i$ parts equal to $k-i$, for $i=0,\dots,k-1$,  
is called an irreducible partition.  Otherwise, the partition 
is reducible.

\begin{definition}
A $k$-Schur function indexed by an irreducible partition 
is said to be a irreducible.
Otherwise, the $k$-Schur function is called reducible.
\label{defindecom}
\end{definition}

\noindent For example, the irreducible $k$-Schur functions for $k=1,2,3$ are
\begin{equation}
\begin{split}
k=1 \, &:  \qquad s_0^{(1)} \, ,\\
k=2 \, &: \qquad s_0^{(2)} \, , \quad s_1^{(2)} \, ,\\
k=3 \, &: \qquad s_0^{(3)} \, , \quad s_1^{(3)}\, ,
\quad s_{2}^{(3)}\, , \quad s_{1,1}^{(3)} \, , \quad s_{2,1}^{(3)} 
\, , \quad s_{2,1,1}^{(3)}\, .
\end{split} \end{equation}

\noindent These examples suggest the following property;

\begin{property} \cite{[LLM]}  
There are $k!$ distinct $k$-irreducible partitions.
\label{lemmakfac}
\end{property}
The concept of irreducibility arose in our study of
the multiplication of a $k$-Schur function with a Schur 
function indexed by a $k$-rectangle.  
It happens that this produces a single $k$-Schur function.
More precisely, it was shown that

\medskip

\begin{theorem}
\cite{[LM2]} 
If $\lambda$ is a $k$-bounded partition, then
\begin{equation}
s_{(\ell^{k-\ell+1})}[X] \, s_{\lambda}^{(k)}[X] = 
s_{\lambda\cup(\ell^{k-\ell+1})}^{(k)}[X] \, .
\end{equation}
\label{conjrecschur}
\end{theorem} 
\noindent Therefore, given the $k!$ irreducible 
$k$-Schur functions,  any reducible $s_\lambda^{(k)}$
is obtained simply by the multiplication of a sequence 
of Schur functions indexed by $k$-rectangles on the proper irreducible
$k$-Schur function.
The case that the $k$-rectangle is the partition $(k)$
relies only on properties of the $k$-split polynomials.
We include the proof here.

\medskip

\begin{property}
For $\lambda$ any $k$-bounded partition,
\begin{equation}
s_k[X]\,s_{\lambda}^{(k)}[X] = 
s_{k,\lambda}^{(k)}[X]\, . 
\end{equation}
\end{property}
\noindent{\bf Proof.}\quad
The unitriangular relation 
between the Schur functions and
$G_\lambda^{(k)}[X]$ given in \eqref{unisG1proof} implies
\begin{equation}
s_k[X]\,s_\lambda^{(k)}[X] = s_k\left(G_\lambda^{(k)}[X] + 
\sum_{\mu>\lambda;\mu_1=\lambda_1} v_{\mu \lambda}^{(k)}\,
G_\mu^{(k)}[X]\right)\,.
\end{equation}
Property~\ref{kpol3} then gives the action of
$s_k[X]$ on a $k$-split polynomial, and we have
\begin{eqnarray}
\label{stup}
s_k[X]\, s_{\lambda}^{(k)}[X] 
& = &
G_{k,\lambda}^{(k)}[X] + 
\sum_{\mu>\lambda;\mu_1=\lambda_1} v_{\mu \lambda}^{(k)} \,
G_{k,\mu}^{(k)}[X]
\, .
\end{eqnarray}
Since each polynomial in the right hand side
of this expression is indexed by a partition
with first component $k$, the expression is
invariant under the action of the
projection operator $T_k^{(k)}$.
That is,
$$
T_k^{(k)}
s_k[X]\, s_{\lambda}^{(k)}[X] 
= 
T_k^{(k)}\left(
G_{k,\lambda}^{(k)}[X] + 
\sum_{\mu>\lambda;\mu_1=\lambda_1} v_{\mu \lambda}^{(k)} \,
G_{k,\mu}^{(k)}[X]
\right)
=
G_{k,\lambda}^{(k)}[X] + 
\sum_{\mu>\lambda;\mu_1=\lambda_1} v_{\mu \lambda}^{(k)} \,
G_{k,\mu}^{(k)}[X]
$$
The right hand side in this expression is the same as 
that of \eqref{stup} and we thus  have that by definition
$T_{k}^{(k)}s_k[X] s_{\lambda}^{(k)}[X] 
= s_k[X] s_{\lambda}^{(k)}[X]
=s_{k,\lambda}^{(k)}$.
\hfill$\square$

\medskip
 
The role of Schur functions indexed by $k$-rectangles
in the subring $\Lambda^{(k)}$ leads naturally
to the study of  the quotient ring $\Lambda^{(k)}/{\mathcal I}_k$, 
where ${\mathcal I}_k$ denotes the ideal generated by 
$s_{(\ell^{k+1-\ell})}[X]$. 
It is known that

\begin{proposition}  
\cite{[LLM]} 
The homogeneous functions indexed by $k$-irreducible partitions
form a basis of the quotient ring $\Lambda^{(k)}/{\mathcal I}_k$. 
\end{proposition}
\noindent Thus, the
dimension of the quotient ring $\Lambda^{(k)}/{\mathcal I}_k$ is $k!$.  
Since we have shown that the $k$-Schur functions form a basis for
$\Lambda^{(k)}$, Theorem~\ref{conjrecschur} implies

\medskip

\begin{theorem}
The irreducible $k$-Schur functions form a basis of the quotient ring 
$\Lambda^{(k)}/{\mathcal I}_k$.
\end{theorem}
In fact, the irreducible $k$-Schur function basis offers
a very beautiful way to carry out 
operations in this quotient ring:  first work in
$\Lambda^{(k)}$ using $k$-Schur functions and then
replace by zero all the $k$-Schur functions indexed by partitions which are
not $k$-irreducible.

\medskip

\section{Analogs of Schur function properties} \label{sec5}

Computer experimentation reveals that many of the properties 
making the Schur function basis so important are generalized
by the $k$-Schur functions.  We now state several of these properties.

\subsection{The $k$-conjugation of a partition}
We give a generalization of partition conjugation
that is an involution on $k$-bounded partitions,
and reduces to usual conjugation of partitions
for large $k$.

A skew diagram $D$ has hook-lengths bounded by $k$ 
if the hook-length of any cell in $D$ is not larger than $k$. 
For a positive integer $m\leq k$, the
$k$-multiplication $m \times^{(k)} D$
is the skew diagram $\overline D$ obtained by prepending
a column of length $m$ to $D$ such that the number of 
rows of $\overline D$ is as small as possible while ensuring 
that its hook-lengths are bounded by $k$.  For example, 
\begin{equation}
{\tiny{\tableau*[scY]{\tf \cr\tf \cr\tf \cr\tf \cr}}} 
\, \times^{(5)} \,  
{\tiny{\tableau*[scY]{  \cr  \cr  & \cr \bl & 
\cr \bl & & & \cr \bl & \bl & \bl & & }}} \, = \,  
{\tiny{\tableau*[scY]{\tf \cr\tf \cr\tf & \cr\tf & \cr\bl & & \cr\bl & \bl & \cr \bl & \bl & & & \cr\bl & \bl & \bl & \bl & & }}} \, .
\end{equation}

\begin{definition}
The $k$-conjugate of a $k$-bounded partition 
$\lambda=(\lambda_1,\dots,\lambda_n)$,
denoted $\lambda^{\omega_k}$, is the vector
obtained by reading the number of boxes in 
each row of the skew diagram,
\begin{equation}
D =  \lambda_1 \times^{(k)} \cdots  \times^{(k)}
\lambda_n  \, ,
\label{spmul}
\end{equation}
arising by $k$-multiplying the entries of $\lambda$ from right to left.
\label{defconju}
\end{definition}
When $k \to \infty$, $\lambda^{\omega_k}=\lambda'=D$ 
since each $k$-multiplication step reduces to 
concatenating a column of height $\lambda_i$.
Further, the $k$-conjugate is an involution on $k$-bounded
partitions:
\begin{theorem}  
\cite{[LLM]}
$\omega_k$ is an involution on partitions bounded by $k$.  
That is, for $\lambda$ with $\lambda_1\leq k$,
\begin{equation}
\left( \lambda^{\omega_k} \right)^{\omega_k} = \lambda \, .
\end{equation}
\end{theorem}

We have observed that the $k$-conjugation of a partition 
plays a natural role in the generalization of 
classical Schur function properties.  We now
give two examples.

\subsection{The involution $\omega$}
It is known that the involution in \eqref{homomor}
acts on a Schur function by 
\begin{equation}
\label{alr}
\omega s_\lambda[X]=s_{\lambda'}[X]\,.
\end{equation}
We thus naturally examine the action of $\omega$ 
on a $k$-Schur function since 
$\omega$ preserves the space $\Lambda^{(k)}$ by
\eqref{ringhk}.
Many examples support the following natural generalization of \eqref{alr}:

\medskip

\begin{conjecture}\label{conjuschur}
 For any $k$-bounded partition $\lambda$,
\begin{equation}
\omega \, s_{\lambda}^{(k)}[X] =  s_{\lambda^{\omega_k}}^{(k)}[X] \, .
\end{equation}
\end{conjecture}
\noindent 
The conjecture holds when $h_M(\lambda)\leq k$ since, 
in this case, $\lambda^{\omega_k}=\lambda'$
and by Property~\ref{largek}, $s_\lambda^{(k)}=s_\lambda$.

\subsection{Pieri Rules}
Beautiful combinatorial algorithms are known for the
Littlewood-Richardson coefficients that
appear in a product of Schur functions;
\begin{equation}
s_{\lambda}[X] \, s_{\mu}[X] = 
\sum_{\nu} c_{\lambda \mu}^{\nu} \, s_{\nu}[X] \,
\quad\text{where}\quad
c_{\lambda\mu}^\nu\in \mathbb N\, . 
\label{schurpieri}
\end{equation}
Since $\Lambda^{(k)}$ is a ring, 
and we have shown that $s_\lambda^{(k)}$ forms a basis 
for this space, a similar expression holds
for the product of two $k$-Schur functions.  That is, 
for $k$-bounded partitions $\lambda$ and $\mu$,
\begin{equation}
\label{klitric}
s_{\lambda}^{(k)}[X] \, s_{\mu}^{(k)}[X] =
\sum_{\nu} {c_{\lambda \mu}^{\nu}}^{\!\!(k)} s_{\nu}^{(k)}[X]  
\quad\text{where}\quad
{c_{\lambda \mu}^{\nu}}^{\!\!(k)} \in \mathbb Z\, , 
\end{equation}
where the integrality of ${c_{\lambda\mu}^{\nu}}^{\!\!(k)}$
follows from Property~\ref{unish1}.
Further, Property~\ref{largek} says that the $k$-Schur functions 
are simply the Schur functions 
when $k$ is large, and therefore
${c_{\lambda \mu}^{\nu}}^{\!\!(k)} = {c_{\lambda \mu}^{\nu}}$
for $k \geq |\nu|$.  In fact,
we believe the coefficients are nonnegative for all $k$.  That is,

\medskip

\begin{conjecture}  For all $k$-bounded partitions $\lambda,\mu,\nu$,
we have $0\leq{c_{\lambda \mu}^{\nu}}^{\!\!(k)} \leq c_{\lambda \mu}^{\nu}$ .
\end{conjecture}

In particular, \eqref{klitric} reduces to a 
$k$-generalization of the Pieri rule
when $\lambda$ is a row (resp. column) of  
length $\ell \leq k$
since $s_{\lambda}^{(k)}[X]$ reduces to $h_{\ell}[X]$
(resp. $e_\ell[X]$).  That is, for $\ell \leq k$,
\begin{equation}
\label{natex}
h_\ell[X]\, s_{\lambda}^{(k)}[X] = \sum_{\mu \in E_{\lambda,\ell}^{(k)}}
s_{\mu}^{(k)}[X]
\quad\text{and}
\quad
e_\ell[X]\, s_{\lambda}^{(k)}[X] = \sum_{\mu \in \bar E_{\lambda,\ell}^{(k)}}
s_{\mu}^{(k)}[X] \, ,
\end{equation}
for some sets of partitions 
$E_{\lambda,\ell}^{(k)}$ and $\bar E_{\lambda,\ell}^{(k)}$,
which we believe naturally extend the Pieri rules by:
\begin{conjecture}
For any positive integer $\ell \leq k$,
\begin{equation}
\begin{split}
E_{\lambda,\ell}^{(k)}\, & = \,
\left\{\mu \, | \,
\mu/\lambda\; is \; a \; horizontal \;\ell\text{-}strip
\;\; and \;\;
\mu^{\omega_k}/\lambda^{\omega_k} \; is \; a \;  vertical
\; \ell\text{-}strip \right\} \, , \\
\bar E_{\lambda,\ell}^{(k)}\, & = \,
\left\{\mu \, | \, 
\mu/\lambda\; is \; a \; vertical \;\ell\text{-}strip
\;\; and \;\;
\mu^{\omega_k}/\lambda^{\omega_k} \; is \; a \;  horizontal
\; \ell\text{-}strip \right\} \, .
\end{split}
\end{equation}
\label{conjpieri}
\end{conjecture}
\noindent 
For example, to obtain the indices of elements
occurring in $e_2 \, s_{3,2,1}^{(4)}$,
we find $(3,2,1)^{\omega_4}=(2,2,1,1)$ by definition.
Adding a horizontal 2-strip to (2,2,1,1) in all ways, 
we obtain
(2,2,2,1,1),(3,2,1,1,1),(3,2,2,1) and (4,2,1,1)
of which all are 4-bounded.
Our set then consists of all the $4$-conjugates of these partitions
that leave a vertical 2-strip when 
$(3,2,1)$ is extracted from them. The 4-conjugates are
\begin{equation}
(2,2,2,1,1)^{\omega_4}   =  \,
{\tiny{\tableau*[scY]{ & \cr & & \cr & & \cr }}} \, ,
\,
(3,2,1,1,1)^{\omega_4}
= \, {\tiny{\tableau*[scY]{ \cr & \cr & \cr & & \cr }}}\, ,
\,
(3,2,2,1)^{\omega_4} = \, {\tiny{\tableau*[scY]{ \cr \cr \cr & \cr & & \cr }}} \, ,
\,
(4,2,1,1)^{\omega_4}  =  \, {\tiny{\tableau*[scY]{ \cr \cr \cr \cr \cr & & \cr }}} \,,
\end{equation}
and of these, the first three are such that a
vertical 2-strip remains when $(3,2,1)$ is extracted.  Thus
\begin{equation}
e_2[X] \, s^{(4)}_{3,2,1}[X] \, = \,   s^{(4)}_{3,3,2}[X] +
s^{(4)}_{3,2,2,1}[X]  +  s^{(4)}_{3,2,1,1,1}[X] \, .
\end{equation}

\section{$t$-generalization} \label{sec6}

The ring of symmetric polynomials over rational functions in an 
extra parameter $t$ has proven  to be of interest in many fields of 
mathematics and physics.  
One natural basis of this space is given by the Hall-Littlewood polynomials,
$H_{\lambda}[X;t]$, which provide $t$-analogs of the
homogeneous symmetric functions $h_{\lambda}[X]$.  
Our approach employs vertex operators that
arise in the recursive construction for the 
Hall-Littlewood polynomials \cite{[Ji]}.
These operators can be defined \cite{[ZS]} for $\ell \in \mathbb Z$, by 
\begin{equation} \label{vertexop}
B_\ell = \sum_{i=0}^{\infty} s_{i+\ell}[X] \, s_{i}[X(t-1)]^{\perp} \, , 
\end{equation}
where for $f,g$ and $h$ arbitrary symmetric functions,  
$f^{\perp}$ is such that on the scalar product \eqref{schurscalar},
\begin{equation}
\langle f^{\perp} g,h \rangle = \langle g,fh \rangle \, .
\end{equation}
The operators add an entry to the Hall-Littlewood polynomials, that is,
\begin{equation} 
\label{hallbb}
H_\lambda[X;t] = B_{\lambda_1} H_{\lambda_2,\ldots,\lambda_\ell}[X;t]
\, , \qquad \text{for} \qquad \lambda_1 \geq \lambda_2 \, .
\end{equation} 
Further, since they satisfy the relation
\begin{equation}
\label{commuB}
B_m B_n = t B_n B_m +t B_{m+1} B_{n-1}- B_{n-1} B_{m+1} \, ,
\qquad m,n \in \mathbb Z \, , 
\end{equation}
and $B_{\ell} \cdot 1=0$ if $\ell<0$, their action on 
$\Lambda$ can be computed algebraically.

We now consider a generalization of the subspace 
 $\Lambda^{(k)}$,
given by 
\begin{equation}
\Lambda^{(k)}_t = {\mathcal L} 
\left\{ H_\lambda[X;t] \right\}_{\lambda;\lambda_1\leq k} \, .
\end{equation}
It is clear that $\Lambda_t^{(1)}\subseteq \Lambda_t^{(2)}\subseteq
\cdots \subseteq\Lambda_t^{(\infty)}=\Lambda$ and thus
that these subspaces provide a filtration for $\Lambda$.
Note that $\Lambda^{(k)}_t$ can be equivalently defined as
\begin{equation}
\Lambda^{(k)}_t = {\mathcal L} \left\{ s_\lambda[X/(1-t)] \right\}_{\lambda;\lambda_1\leq k} \, , \qquad \text{or} \quad \Lambda^{(k)}_t = {\mathcal L} \left\{ H_\lambda[X;q,t] \right\}_{\lambda;\lambda_1\leq k} \, .
\end{equation}
Again, we want elements that play the role in $\Lambda^{(k)}_t$
that the Schur functions play in $\Lambda$.  In particular, 
since $\Lambda^{(\infty)}_t=\Lambda$,
we want a basis for $\Lambda_t^{(k)}$ that reduces to 
the Schur functions when $k$ is large.

Our construction of these elements 
will naturally extend Definition~\ref{kschurdef} 
for the $k$-Schur functions in the case $t=1$.  First, multiplication by 
$s_{\lambda_1}[X]$ is replaced by the action of the 
operator $B_{\lambda_1}$.  Then, the projection 
operator $T_{\lambda_1}^{(k)}$ is replaced by a $t$-analogous
operator which, to define, requires an appropriate extension of
the $k$-split polynomials.  Recent developments in the theory of 
symmetric functions aid us with this task.

\subsection{$k$-split polynomials}
Important in our work with the $k$-Schur functions
is a family of polynomials, studied in many recent papers 
such as \cite{[LLT],[Shi],[S2],[S1],[ZS]}, 
that give a $t$-analog of the product of Schur functions.
These functions, indexed by a sequence of partitions, can
be built recursively using vertex operators \cite{[ZS]}.
For a partition $\lambda$ of length $m$, 
define 
\begin{equation}
\label{jtt}
B_\lambda \equiv \prod_{1\leq i<j\leq m}
(1-te_{ij})
B_{\lambda_1} \cdots
B_{\lambda_m} \, ,
\end{equation}
where $e_{ij}$ acts by
\begin{equation}
e_{ij}\left(
B_{\lambda_1} \cdots
 B_{\lambda_m}
\right)
=
B_{\lambda_1} \cdots
B_{\lambda_i+1}
\cdots
B_{\lambda_j-1} \cdots
B_{\lambda_m}
\, .
\end{equation}

For any sequence of partitions $(\lambda^{(1)},\lambda^{(2)},\ldots)$, 
the generalized Schur function product is then defined recursively by
\begin{equation} \label{recursiB}
{\H}_{(\lambda^{(1)},\lambda^{(2)},\lambda^{(3)},\dots)}[X;t] = 
B_{\lambda^{(1)}} {\H}_{(\lambda^{(2)},\lambda^{(3)},\dots)} [X;t] \, ,
\end{equation}
starting with ${\H}_{()}=1$.
Note that since $B_{\lambda} \cdot 1= s_{\lambda}[X]$, we have 
that
\begin{equation}
\label{Bon1}
{\H}_{(\lambda)}[X;t]=s_{\lambda}[X] \, .
\end{equation}
Appendix \ref{A1} gives an earlier formulation \cite{[S1]} for 
${\H}_{(\lambda^{(1)},\lambda^{(2)},\dots)}[X;t]$.  
Note in \cite{[ZS]}, $B_{\lambda}$ is denoted $H_{\lambda}^t$ and is
given in terms of generating series.  Formula \eqref{jtt} can be extracted
from their formula (17).

It was shown \cite{[S1]} that for a sequence of partitions,
$S=( \lambda^{(1)},\lambda^{(2)},\cdots)$,
\begin{equation} 
{\H}_{S}[X;t] = 
\sum_{\mu\vdash|\lambda^{(1)}|+|\lambda^{(2)}|+\cdots}K_{\mu;S}(t) s_\mu[X] \, ,
\quad\text{where}\quad K_{\mu;S}(t)\in\mathbb Z[t]\, .
\end{equation}
Since the action of $B_{\lambda}$ on an arbitrary function $f$
reduces to multiplication by $s_{\lambda}[X]$
\begin{equation}\label{limt1}
B_{\lambda}\, f = s_{\lambda}[X] \, f \, , \qquad \text{when~} t=1 \, ,
\end{equation}
the $K_{\mu;S}(t)$ are known as generalized Kostka 
polynomials since, in this case, they satisfy
\begin{equation}
K_{\mu;S}(1)= \langle s_{\mu}[X], s_{\lambda^{(1)}}[X] s_{\lambda^{(2)}}[X] 
\cdots \rangle \, .
\end{equation}

For our purposes, we consider only the ${\H}_S[X;t]$ indexed by a 
dominant sequence $S$.  That is, sequences of partitions 
$S=(\lambda^{(1)},\lambda^{(2)},\ldots)$
such that the concatenation of $\lambda^{(1)},\lambda^{(2)},\dots$, 
denoted $\bar S$, forms a partition.  In this case, we 
prove that ${\H}_S[X;t]$ obeys important unitriangular relations.

\medskip

\begin{property} \label{propHtri}
If $S$ is dominant, with $\bar S=\lambda$, then
\begin{equation}
{\H}_S[X;t]= s_{\lambda}[X] + \sum_{\mu > \lambda} K_{\mu;S}(t) s_{\mu}[X] \, ,
\quad\text{where}\quad K_{\mu;S}(t)\in\mathbb Z[t] \, ,
\label{Hins}
\end{equation}
\begin{equation}
{\H}_S[X;t]= H_{\lambda}[X] + \sum_{\mu > \lambda} C_{\mu;S}(t) H_{\mu}[X] \, ,
\label{HinHall}
\quad\text{where}\quad C_{\mu;S}(t)\in\mathbb Z[t]
\, .
\end{equation}
\end{property}
\noindent{\bf Proof.}\quad
Proposition~\ref{propuni}, with Proposition~\ref{propsmaller},
prove relation \eqref{Hins} (see Appendix~\ref{A2}).
The second identity then follows from unitriangularity \eqref{HallinS} 
of $s_\lambda[X]$ in terms of $H_\mu[X;t]$.  \hfill$\square$

\medskip

We have discovered that a particular subset of the 
$\mathcal H_{S}$ not only form a basis for $\Lambda^{(k)}_t$,
but are essential in 
the construction of the  $k$-Schur functions.

\medskip

\begin{definition}
\label{Defksp}
The $k$-split polynomials are defined,
for a $k$-bounded partition $\lambda$, by
\begin{equation}
G_\lambda^{(k)}[X;t] = {\H}_S[X;t] 
\quad\text{where}
\quad S=\lambda^{\to k}\quad\text{is the $k$-split of $\lambda$}\,. 
\label{defksp}
\end{equation}
\end{definition}

Since ${\H}_S[X;t]$ reduces to the product of Schur
functions indexed by elements of $S$ when $t=1$, 
\begin{equation} \label{limt1G}
G_{\lambda}^{(k)}[X;1]=G_{\lambda}^{(k)}[X] \, ,  
\end{equation}
and therefore, 
$G_\lambda^{(k)}[X;t]$ is a proper $t$-generalization of 
the $G_\lambda^{(k)}[X]$ introduced in Section~\ref{sec3}.

\medskip

It develops that the $k$-split polynomials satisfy 
properties analogous to those held by the $k$-split polynomials at $t=1$.
We first show that the $G_\lambda^{(k)}[X;t]$ actually
lie in the space $\Lambda_t^{(k)}$.  To this end, we start by 
showing that operators $B_i$ preserve this space. 

\begin{proposition} \label{coropreserve}
If $f \in \Lambda^{(k)}_t$ then
$B_i \, f \in \Lambda^{(k)}_t$
for all $i\in\mathbb Z$ with $i \leq k$.
\end{proposition}

This claim follows from a preliminary result
on subspaces of $\Lambda^{(k)}_t$.  More precisely,
for $a \in \mathbb Z$ and $b \in \mathbb N$, with $a \leq  b$, 
we define 
\begin{equation}
\Lambda^{(a,b)}_t = {\mathcal L}\{ H_\lambda[X;t] \}_{a\leq \lambda_1\leq b}
\,, 
\end{equation}
and $\Lambda^{(a,b)}_t\!=\!\Lambda^{(0,b)}_t$ for $a<0$.
In particular, $\Lambda^{(k)}_t\!=\!\Lambda^{(0,k)}_t$ and
$\Lambda^{(a,b)}_t \subseteq \Lambda^{(k)}_t$ if $b \leq k$.
Thus, Proposition~\ref{coropreserve} is an immediate consequence of

\begin{lemma}
\label{lemplusgrand}
If $f \in \Lambda^{(k)}_t$ then $B_i \, f \in \Lambda^{(i,k)}_t$
for all integers $i \leq k$.
\end{lemma}
\noindent{\bf Proof.}\quad
The assertion holds for $i=k$ since
$B_k H_{\lambda}[X;t]=H_{(k,\lambda_1,\lambda_2,\dots)}[X;t]\in
\Lambda^{(k,k)}_t$.  Assume by induction that for 
all $f\in \Lambda_t^{(k)}$, 
\begin{equation}
\label{ind1}
B_i \, f \in \Lambda^{(i,k)}_t
\qquad j < i \leq k\, .
\end{equation}
Thus, it suffices to show that $B_j H_\lambda[X;t]\in\Lambda^{(j,k)}_t$.

We prove this claim by induction on $\ell(\lambda)$.
First, $B_j H_{()}[X;t]=H_{(j)}[X;t]\in\Lambda^{(j,k)}_t$ 
for $j\geq 0$ and otherwise $B_j H_{()}[X;t]=0 \in \Lambda^{(0,k)}_t 
= \Lambda^{(j,k)}$.
Now, let $\lambda=(\lambda_1,\ldots,\lambda_n)$ be a $k$-bounded partition
and assume $B_j H_{\hat\lambda}[X;t]\in\Lambda^{(j,k)}_t$
for all $\hat\lambda=(\lambda_2,\ldots,\lambda_n)$.
If $\lambda_1\leq j$ then
$B_j H_{\lambda}[X;t]=H_{(j,\lambda_1,\lambda_2,\dots)}
\in\Lambda^{(j,k)}_t$ by \eqref{hallbb}.
Now consider $\lambda_1\!=\!j+1$. 
The commutation relation \eqref{commuB}
reduces to $B_j H_{\lambda}=B_j B_{j+1} H_{\hat \lambda} 
= t B_{j+1} B_j H_{\hat\lambda}$.  
Since $f=B_j H_{\hat\lambda}\in\Lambda^{(j,k)}_t\subseteq\Lambda_t^{(k)}$ 
by assumption, $B_j H_{\lambda}
=t B_{j+1} f\in\Lambda^{(j+1,k)}_t\subseteq\Lambda_t^{(j,k)}$
by \eqref{ind1}.  Finally, for $\lambda_1>j+1$ we again use
\eqref{commuB} to obtain
$B_jH_{\lambda} = B_j B_{\lambda_1} H_{\hat \lambda}= 
t B_{\lambda_1} B_j H_{\hat\lambda} + 
t B_{j+1} B_{\lambda_1-1} H_{\hat\lambda}
- B_{\lambda_1-1} B_{j+1} H_{\hat\lambda}$.  
Regarding the first term in the right hand side,
our assumption gives that
$B_j H_{\hat\lambda}\in\Lambda^{(j,k)}_t\subseteq\Lambda_t^{(k)}$ 
and then using \eqref{ind1},
$tB_{\lambda_1} B_j H_{\hat\lambda}\in\Lambda_t^{(\lambda_1,k)}
\subseteq\Lambda_t^{(j,k)}$ since $\lambda_1>j$.
Similar reasoning applies to the second and third term,
and we thus have $B_j H_{\lambda}\in \Lambda^{(j,k)}_t$,
proving our claim.  \hfill $\endprf$

Now, given that $\Lambda_t^{(k)}$ is invariant under $B_i$,
we can prove a more general statement that will
imply the $k$-split polynomials lie in $\Lambda^{(k)}_t$.

\medskip

\begin{proposition}
\label{lempreserve}
If $\lambda$ is a partition with $h_M(\lambda)\leq k$,
then $B_\lambda\, f\in \Lambda_t^{(k)}$ 
for any
$f\in\Lambda_t^{(k)}$. 
\end{proposition}
\noindent {\bf{Proof.}} 
An operator version of the Jacobi-Trudi determinant \eqref{jt}
is given by \eqref{jtt}. 
Thus, the expansion of $B_\lambda$ in terms of products of 
the operators $B_i$ yields only terms $B_i$ with $i\leq k$.
Proposition~\ref{coropreserve} then implies the result.
\hfill $\endprf$

\medskip

\begin{property}  \label{GinLambda}
For any $k$-bounded partition $\lambda$, we have that
\begin{equation}
G_{\lambda}^{(k)}[X;t] \in \Lambda^{(k)}_t
\, .
\end{equation}
\end{property}
\noindent{\bf Proof.}\quad
Recursion \eqref{recursiB} for $\H_S$ allows us to 
restate Definition~\ref{Defksp} as,
\begin{equation}
\label{defksp2}
G_\lambda^{(k)}[X;t] = 
B_{\lambda^{(1)}} 
\cdots
B_{\lambda^{(n)}} 
\cdot 1 \, ,
\quad\text{where}\quad
\lambda^{\to k} = (\lambda^{(1)},\cdots,\lambda^{(n)}).
\end{equation}
Since all $h_M(\lambda^{(i)})\leq k$,
Proposition~\ref{lempreserve} implies
$B_{\lambda^{(n)}}\cdot 1 \in\Lambda_t^{(k)}$
given $1 \in \Lambda^{(k)}_t$. 
By induction,
assuming $B_{\lambda^{(2)}}\cdots B_{\lambda^{(n)}}\cdot 1\in\Lambda_t^{(k)}$,
Proposition~\ref{lempreserve} verifies our claim.
\hfill$\square$

\medskip

\begin{property}  \label{propGH}
We have
\begin{equation}
G_{\lambda}^{(k)}[X;t] = H_{\lambda}[X;t] + 
\sum_{\mu > \lambda; \mu_1 \leq k} g_{\lambda \mu}^{(k)}(t) \, 
H_{\mu}[X;t] \, , \quad
\text{where} \quad g_{\lambda \mu}^{(k)}(t) \in \mathbb Z[t]\, .
\end{equation}
\end{property}
\noindent{\bf Proof.}\quad
The $k$-split polynomials can be expanded in terms of
Hall-Littlewood polynomials indexed by $k$-bounded partitions
since Property~\ref{GinLambda} proves they lie in $\Lambda_t^{(k)}$.
Then, by the unitriangular expansion of $\H_S$ given in \eqref{HinHall},
$G_\lambda^{(k)}[X;t]=\H_{\lambda^{\to k}}[X;t]$ has the asserted 
unitriangularity. 
\hfill$\square$

\medskip

\begin{property} \label{propsplitschur}
 Let $\lambda$ be such that $h_M(\lambda) \leq k$.  Then,
\begin{equation}
G_{\lambda}^{(k)}[X;t] = s_{\lambda}[X] \, .
\end{equation}
\end{property}
\noindent{\bf Proof.}\quad 
$G_{\lambda}^{(k)}[X;t]=\H_{(\lambda)}[X;t]$
since $\lambda^{\to k} = (\lambda)$
for $h_M(\lambda) \leq k$.
The claim then follows from \eqref{Bon1}, 
which states that $\H_{(\lambda)}[X;t]=s_{\lambda}[X]$
for all $\lambda$.
\hfill $\endprf$

\medskip

\begin{theorem}
The $k$-split polynomials form a basis
of $\Lambda_t^{(k)}$.
\end{theorem}
\noindent{\bf Proof.}\quad 
Since $\Lambda^{(k)}_t$ is the span of Hall-Littlewood
polynomials indexed by $k$-bounded partitions
and $G_\lambda^{(k)}$ are also indexed by $k$-bounded partitions,
the theorem follows from Property~\ref{propGH}.
\hfill$\square$

\medskip

\subsection{$k$-Schur functions}  
As with the $t=1$ case, although the $k$-split polynomials 
form a basis for $\Lambda^{(k)}_t$, they
do not play the fundamental role that the Schur functions do 
for $\Lambda$.  However, these polynomials are needed in
the construction of our Schur analog, $s_\lambda^{(k)}[X;t]$.  
We use a projection operator, ${\bar T}_j^{(k)}$, for $j \leq k$, that acts linearly on $\Lambda^{(k)}_t$ by
\begin{equation}
{\bar T}_j^{(k)} \, G_\lambda^{(k)}[X;t] = 
\begin{cases}
G_\lambda^{(k)}[X;t] & \text{if } \lambda_1 = j\\
0 & \text{otherwise }
\end{cases}
\, .
\end{equation}

\medskip

\begin{definition}  
\label{defks}
For $k$-bounded partition $\lambda$,
the $k$-Schur functions are recursively defined
\begin{equation}
s_{\lambda}^{(k)}[X;t] = {\bar T}_{\lambda_1}^{(k)} B_{\lambda_1} 
s_{(\lambda_1,\lambda_2,\dots)}^{(k)}[X;t] \, ,
\quad\text{where}\quad
s_{()}^{(k)}[X;t]=1\, .
\end{equation}
\end{definition}
\noindent
Tables of $k$-Schur functions in terms of Schur functions
can be found in Subsection~\ref{tabkschur}.

When $t=1$, from (\ref{limt1}) and (\ref{limt1G}), 
we recover the $k$-Schur functions of Definition~\ref{kschurdef}, that is,
\begin{equation}
s_\lambda^{(k)}[X;1] = s_\lambda^{(k)}[X] \, .
\end{equation}
We now prove several other properties satisfied by the
 $k$-Schur functions.

\medskip

\begin{property} 
\label{propstri}
For $\lambda$ a $k$-bounded partition, we have
\begin{equation}
s_{\lambda}^{(k)}[X;t] = G_{\lambda}^{(k)}[X;t] + 
\sum_{\mu > \lambda\atop\mu_1=\lambda_1} u_{\mu \lambda}^{(k)}(t) 
\, G_{\mu}^{(k)}[X;t] \, ,
\quad\text{where}\quad u_{\mu \lambda}^{(k)}(t) \in \mathbb Z[t]\, .
\end{equation}
\end{property}
\noindent{\bf Proof.}\quad 
The assertion holds for $s_{()}^{(k)}[X;t]=1=G_{()}^{(k)}[X;t]$ 
by definition, Let $\lambda=(\lambda_1,\dots,\lambda_n)$ be 
a $k$-bounded partition and assume
by induction on $n$ that the property is 
true for $\hat \lambda=(\lambda_2,\dots,\lambda_n)$.  
Definition~\ref{defks} gives that
$s_{\lambda}^{(k)}[X;t] = {\bar T}_{\lambda_1}^{(k)} B_{\lambda_1} 
s_{\hat \lambda}^{(k)}[X;t]$.  By the induction hypothesis,
we have
\begin{equation}
\begin{split}
B_{\lambda_1} 
s_{\hat \lambda}^{(k)}[X;t]
&= 
 B_{\lambda_1} \left( G_{\hat \lambda}^{(k)}[X;t] + 
\sum_{\mu >\hat  \lambda; \mu_1=\hat \lambda_1} u_{\mu \hat \lambda}(t) 
\, G_{\mu}^{(k)}[X;t] \right) \, , 
\qquad u_{\mu\hat\lambda}(t)\in\mathbb Z[t]\, .
\end{split}
\end{equation}
Recall the $k$-split polynomial is defined 
\eqref{defksp} by $G_\gamma^{(k)}=\H_{\gamma^{\to k}}$,
and thus we have
\begin{equation}
\begin{split}
B_{\lambda_1} s_{\hat \lambda}^{(k)}[X;t]
& = B_{\lambda_1} \left( \H_{\hat \lambda^{\to k}}[X;t] + 
\sum_{\mu >\hat  \lambda; \mu_1=\hat \lambda_1} u_{\mu \hat \lambda}(t) 
\, \H_{\mu^{\to k}}[X;t] \right) \\
& = 
 \H_{\bigl((\lambda_1),\hat \lambda^{\to k}\bigr)}[X;t] + 
\sum_{\mu >\hat  \lambda; \mu_1=\hat \lambda_1} u_{\mu \hat \lambda}(t) 
\, \H_{\bigl((\lambda_1),\mu^{\to k}\bigr)}[X;t] \, .
\end{split}
\end{equation}
$\mu_1\!=\!\hat\lambda_1\!=\!\lambda_2\!\leq\!\lambda_1$
implies that $\bigl( (\lambda_1),\hat \lambda^{\to k} \bigr)$
and $\bigl((\lambda_1), \mu^{\to k}\bigr)$ are dominant
sequences.  Further, $\hat \lambda < \mu$
implies that
$\lambda=(\lambda_1) \cup \hat \lambda <
(\lambda_1) \cup  \mu$.
Therefore, unitriangularity of $\H_S$ \eqref{HinHall}
further gives that
\begin{equation}
B_{\lambda_1} s_{\hat \lambda}^{(k)}[X;t]
= 
H_{ \lambda}[X;t] + 
\sum_{\mu > \lambda} v_{\mu  \lambda}(t) 
\, H_{\mu}[X;t]\, , \qquad  v_{\mu  \lambda}(t) \in \mathbb Z[t] 
\, . 
\label{pre}
\end{equation}
$s_{\hat \lambda}^{(k)}[X;t] \in \Lambda^{(k)}_t$ 
since the induction hypothesis gives that it can be expanded 
in terms of $G_\mu^{(k)}[X;t]$.
Thus, $B_{\lambda_1}s_{\hat \lambda}^{(k)}[X;t]\in\Lambda^{(k)}_t$ 
by Proposition~\ref{coropreserve}
and therefore the coefficients $v_{\mu\lambda}(t)$ are
non-zero only for $\mu_1\leq k$. \eqref{pre} then becomes
\begin{equation}
\begin{split}
 B_{\lambda_1} 
s_{\hat \lambda}^{(k)}[X;t] & =  H_{ \lambda}[X;t] + 
\sum_{\mu > \lambda; \mu_1 \leq k} v_{\mu  \lambda}(t) 
\, H_{\mu}[X;t] 
\, , 
\qquad  v_{\mu  \lambda}(t) \in \mathbb Z[t] 
\end{split} 
\end{equation}
Unitriangularity between Hall-Littlewood polynomials
and   $k$-split polynomials gives,
\begin{equation}
\begin{split}
B_{\lambda_1} 
s_{\hat \lambda}^{(k)}[X;t]
& =  G_{ \lambda}^{(k)}[X;t] + 
\sum_{\mu > \lambda; \mu_1 \leq k} u_{\mu  \lambda}(t) 
\, G_{\mu}^{(k)}[X;t] \, , \qquad  u_{\mu  \lambda}(t) \in \mathbb Z[t] \, .
\end{split} 
\end{equation}
Applying ${\bar T}_{\lambda_1}^{(k)}$ to both sides of this expression,
we arrive at
\begin{equation}
\begin{split}
{\bar T}_{\lambda_1}^{(k)}  B_{\lambda_1} 
s_{\hat \lambda}^{(k)}[X;t] 
& =
G_{ \lambda}^{(k)}[X;t] + 
\sum_{\mu > \lambda; \mu_1 = \lambda_1} u_{\mu  \lambda}(t) 
\, G_{\mu}^{(k)}[X;t] \, ,
\qquad u_{\mu  \lambda}(t) \in \mathbb Z[t]\, .
\end{split}
\end{equation}
Since 
$s_{\lambda}^{(k)}= 
\bar T_{\lambda_1}^{(k)} B_{\lambda_1}s_{\hat\lambda}^{(k)}$,
we have our claim.
\hfill $\endprf$

\medskip

The following theorem is an immediate consequence of 
this property.

\medskip

\begin{theorem}  
The $k$-Schur functions form a basis of $\Lambda^{(k)}_t$.
\end{theorem}

We can thus refine the expansion of Hall-Littlewood
polynomials in terms of Schur functions \eqref{HallinS}. 
That is, for any $k$-bounded partition $\lambda$,
\begin{equation}
H_\lambda[X;t] = 
s_\lambda^{(k)}[X;t] + 
\sum_{\mu>\lambda: \mu_1\leq k} K_{\mu\lambda}^{(k)}(t)\,s_\mu^{(k)}[X;t]\, ,
\quad\text{where}\quad K_{\mu\lambda}^{(k)}(t)\in\mathbb Z[t]\, .
\end{equation}
The integrality of $K_{\mu\lambda}^{(k)}(t)$ follows from the
unitriangularity and integrality  in 
Properties~\ref{propGH} and \ref{propstri}.
Moreover, by our triangularity and integrality properties and \eqref{HallinS},
we also have integrality of the coefficients in
Conjecture~\eqref{posischur}:
\begin{property}
For any $k$-bounded partition $\lambda$,
\begin{equation} 
s_{\lambda}^{(k)}[X;t]= s_\lambda[X]
+\sum_{\mu >\lambda} v_{\mu \lambda}^{(k)}(t)\,
s_{\mu}[X]\, ,\quad\text{where}\quad
v_{\mu \lambda}^{(k)}(t) \in \mathbb Z[t] \, .
\end{equation} 
\end{property}
This unitriangularity property, given that the coefficients in the
Schur function expansion of the Macdonald polynomials
are polynomials in $q$ and $t$ with integral coefficients,
implies
\begin{property}
For any $k$-bounded partition $\lambda$,
\begin{equation} 
H_{\lambda}[X;q,t]= 
\sum_{\mu: \mu_1\leq k } 
K_{\mu \lambda}^{(k)}(q,t)\,
s_{\mu}^{(k)}[X;t]\, ,\quad\text{where}\quad
K_{\mu \lambda}^{(k)}(q,t) \in \mathbb Z[q,t] \, .
\end{equation} 
\end{property}

Now, to further support the idea that the $s_\lambda^{(k)}[X;t]$
provide a refinement for Schur function theory, 
we must show that they reduce to the usual $s_\lambda[X]$
when $k\to\infty$.  This result relies on a lemma.

\begin{lemma} \label{lemsepara}
If $\lambda$ is a partition with $h_M(\lambda) \leq k$
then
${\bar T}^{(k)}_{\lambda_1}\, B_{\lambda} \, f =
{\bar T}^{(k)}_{\lambda_1}\, B_{\lambda_1} \, B_{\hat \lambda} \, f $
for all $f \in \Lambda^{(k)}_t$.
\end{lemma}
\noindent{\bf Proof.}\quad
We need to show that, for $f \in \Lambda^{(k)}_t$,
${\bar T}_{\lambda_1}^{(k)}\left( B_\lambda-B_{\lambda_1}B_{\hat\lambda}
\right)\cdot f=0$.
Definition \eqref{jtt} gives
\begin{equation}
\begin{split}
B_\lambda
& =
\prod_{2\leq j \leq \ell(\lambda)} (1-te_{1j}) B_{\lambda_1}
\prod_{2\leq i <j \leq \ell(\lambda)} (1-te_{ij}) B_{\lambda_2} \cdots
B_{\lambda_{n}}
\\
&=
\prod_{2\leq j \leq \ell(\lambda)} (1-te_{1j}) B_{\lambda_1}
B_{\hat\lambda}
\end{split}
\end{equation}
In the expansion of this product,
with the exception of the term $B_{\lambda_1}B_{\hat\lambda}$,
each term contains at least one $e_{1j}$,
which increases the index of $B_{\lambda_1}$.
Therefore, using the argument of Proposition~\ref{lempreserve},
we have
\begin{equation}
B_{\lambda}= B_{\lambda_1} B_{\hat \lambda}+
\sum_{i=1}^{\ell(\lambda)-1}
c_i(t) \, B_{\lambda_1+i} \, O_i \, ,\quad \qquad c_i(t) \in \mathbb Z[t] \, ,
\label{te}
\end{equation}
where $O_i$ is a product of $B_{j}$'s with $j\leq k$.
Proposition~\ref{coropreserve} states that
$B_i\,f\in\Lambda_t^{(k)}$ for $f\in\Lambda_t^{(k)}$
and $i\leq k$ and it is thus clear
that $O_i \cdot f \in \Lambda^{(k)}_t$.
Furthermore, since $\lambda_1+i\leq \lambda_1+\ell(\lambda)-1=h_M(\lambda)\leq k$,
we have $B_{\lambda_1+i} O_i \cdot f \in \lambda^{(\lambda_1+i,k)}_t
\subseteq\Lambda_t^{(\lambda_1+1,k)}$ by Lemma~\ref{lemplusgrand}.
Thus, from expression \eqref{te}, 
\begin{equation}
\left(B_{\lambda}- B_{\lambda_1} B_{\hat \lambda}\right)\cdot f
=
\sum_{i=1}^{\ell(\lambda)-1}
c_i(t) \, B_{\lambda_1+i} \, O_i \cdot f
\in \Lambda_t^{(\lambda_1+1,k)} \, .
\end{equation}
The definition of $\Lambda_t^{(\lambda_1+1,k)}$,
with unitriangularity in Proposition~\ref{propGH},
then gives the expansions
\begin{equation}
\left(B_{\lambda}- B_{\lambda_1} B_{\hat \lambda}\right)\cdot f
= \sum_{\mu; \lambda_1+1 \leq \mu_1 \leq k} b_\mu(t) \,
H_{\mu}[X;t]
= \sum_{\mu; \mu\geq\lambda_1+1} d_\mu(t) \,
G_{\mu}^{(k)}[X;t] \, .
\end{equation}
Finally, acting with $\bar T_{\lambda_1}^{(k)}$, we have
${\bar T}_{\lambda_1}^{(k)}\left(B_{\lambda_1}
B_{\hat \lambda}\cdot f-B_{\lambda}\cdot f \right)= 0$.
\hfill $\square$

\medskip

\begin{property} \label{kgrandt}
For any $k$-bounded partition $\lambda$, 
\begin{equation}
s_{\lambda}^{(k)}[X;t] = s_{\lambda}[X] \, \quad 
\text{ if } \quad h_M(\lambda) \leq k \, .
\end{equation}
\end{property}
\noindent{\bf Proof.}\quad
We proceed by induction on $\ell(\lambda)$.
First, $s_{()}^{(k)}[X;t]=1=s_{()}[X]$. 
Let $\lambda=(\lambda_1,\dots,\lambda_n)$ be 
a $k$-bounded partition with $h_M(\lambda)\leq k$.
If $\hat\lambda=(\lambda_2,\dots,\lambda_n)$  then 
$h_M(\hat\lambda)\leq h_M(\lambda)$ and induction 
implies $s_{\hat\lambda}^{(k)}[X;t]=s_{\hat\lambda}[X]$.
Thus, $s_{\lambda}^{(k)}[X;t]={\bar T}_{\lambda_1}^{(k)}B_{\lambda_1}
s_{\hat\lambda}[X]$ and it  suffices to show 
${\bar T}_{\lambda_1}^{(k)}B_{\lambda_1}s_{\hat\lambda}[X]=s_\lambda[X]$.
Since \eqref{Bon1} gives that $B_{\mu} \cdot 1= s_{\mu}[X]$, 
we need only show 
${\bar T}_{\lambda_1}^{(k)}B_{\lambda_1}B_{\hat\lambda}\cdot 1=
B_\lambda\cdot 1$, or 
${\bar T}_{\lambda_1}^{(k)}\left( B_\lambda-B_{\lambda_1}B_{\hat\lambda}
\right)\cdot 1=0$.  
This result follows from the special case $f=1 \in \Lambda^{(k)}_t$ of the 
previous lemma.
\hfill $\endprf$

\section{Irreducibility and $k$-conjugation} \label{sec7}

\noindent The concept of irreducibility introduced in Section~\ref{secirred} 
extends naturally in the space $\Lambda_t^{(k)}$.  In fact,
the following $t$-analog of Theorem~\ref{conjrecschur}
holds:
\begin{theorem}
\cite{[LM2]}:
For any $k$-bounded partition $\lambda$,
\begin{equation}
B_{(\ell^{k-\ell+1})}\, s_{\lambda}^{(k)}[X;t]
=t^*\,s_{\lambda \cup (\ell^{k-\ell+1})}^{(k)}[X;t] \, ,
\end{equation}
where $t^*$ stands for some power of $t$.
\end{theorem}
This theorem implies that when $k=2$, any $k$-Schur function
can be built using a sequence of operators $B_{2}$ and $B_{1,1}$
applied to either $s_{1}^{(2)}[X;t]=s_1[X]$ or 
$s_{()}^{(2)}[X;t]=s_{()}[X]=1$.
Thus, there is a connection between the $k$-Schur functions 
and the positive functions introduced in \cite{[LM],[Z]} and, 
from \cite{[ZS]}, to the generalized Kostka polynomials 
of \cite{[S1]}.  This connection implies that
our refinement (\ref{posiMac}) of Macdonald's original 
positivity conjecture holds when $k=2$.

With regard to extending Conjecture~\ref{conjuschur} to the
general case, we first note that the involution $\omega$ 
is not well defined on $\Lambda^{(k)}_t$. 
That is, an element $f \in \Lambda^{(k)}_t$ may be such that 
$\omega f \not \in \Lambda^{(k)}_t$.  However, there is a simple 
generalization for $\omega$ that preserves $\Lambda^{(k)}_t$.  
Let $\omega_t$ be defined on an element $f \in \Lambda$, by
\begin{equation}
\omega_t f = (\omega f)\big|_{t \to 1/t} \, .
\end{equation}
For instance, if $f= \sum_{\mu} c_{\mu}(q,t) \, s_{\mu}[X]$, we have
\begin{equation}
\omega_t f = \sum_{\mu} c_{\mu}(q,1/t) \, s_{\mu'}[X] \, .
\end{equation}
Since $\omega$ is an involution, 
$\omega_t$ is also.  In fact, it is well defined
on $\Lambda^{(k)}_t$.
\begin{proposition} If $f \in \Lambda^{(k)}_t$, then 
$\omega_t f \in \Lambda^{(k)}_t$.
\end{proposition}
\noindent 
{\bf Proof.}\quad 
Since $\Lambda^{(k)}_t$ can also be defined as 
$\Lambda^{(k)}_t={\mathcal L}\{s_{\lambda}[X/(1-t)] 
\}_{\lambda;\lambda_1 \leq k}$,
it suffices to show, for any partition $\lambda$,
\begin{equation}
\omega_t \, s_{\lambda}[X/(1-t)]= (-t)^{|\lambda|}s_{\lambda}[X/(1-t)] \, .
\end{equation}
Let ${\bar X}$ denote the alphabet $X$ where $x_i \to -x_i$ for all $i$.  
That is, $P[{\bar X}]=P(-x_1,-x_2,\dots)$
for an arbitrary symmetric function $P[{X}]=P(x_1,x_2,\dots)$. 
Note, this implies that if $P[X]$ is homogeneous of 
degree $d$, then $P[{\bar X}]=(-1)^d P[X]$.  
Since $\omega$ acts on 
$P[X]$ by $\omega P[X] = P[-{\bar X}]$, we thus have
\begin{equation}
\omega_t s_{\lambda}[X/(1-t)]=s_{\lambda}[-{\bar X}/(1-1/t)]=
s_{\lambda}[t {\bar X}/(1-t)]= (-t)^{|\lambda|} s_{\lambda}[X/(1-t)] \, ,
\end{equation}
which completes the proof.  \hfill $\endprf$

The action of $\omega_t$ on a $k$-Schur function appears 
to take a simple form generalizing Conjecture~\ref{conjuschur}.
\begin{conjecture} For any $k$-bounded partition $\lambda$,
\begin{equation}
\omega_t \, s_{\lambda}^{(k)}[X;t] = t^{-c(\lambda)} 
s_{\lambda^{\omega_k}}^{(k)}[X;t] \, ,
\end{equation}
where $c(\lambda)$ is some nonnegative integer.
\end{conjecture}

Another unique property of the Schur functions is the expansion,
\begin{equation}
s_\lambda[X+Y]=\sum_{|\mu|+|\rho|=|\lambda|}
c_{\mu \rho}^\lambda\,s_{\mu}[X]\,s_{\rho}[Y]
\quad\text{where}\quad
c_{\mu\rho}^\lambda\in \mathbb N\, .
\end{equation}
We have found by experimentation that the
$k$-Schur functions also satisfy a similar relation,
\begin{conjecture}
For any $k$-bounded partition,
\begin{equation}
s_\lambda^{(k)}[X+Y;t]=\sum_{|\mu|+|\rho|=|\lambda|}
g_{\mu \rho}^\lambda(t)\,s_{\mu}^{(k)}[X;t]\,s_{\rho}^{(k)}[Y;t]\,
\quad\text{where}
\quad g_{\mu\rho}^\lambda\in \mathbb N[t]\,.
\end{equation}
\end{conjecture}

\section{Appendix} \label{sec8}
\subsection{Definition of $\H_S[X;t]$ \cite{[S1]}} 
\label{A1}
Consider a sequence of partitions
$S=(\lambda^{(1)},\lambda^{(2)},\lambda^{(3)},...)$
with $\eta_i$ parts (some of which may be zero) in $\lambda^{(i)}$.
If $\eta=(\eta_1,\eta_2,\ldots)$ and
$\ell(S) \equiv |\eta|=n$, we set
\begin{equation}
\text{Roots}_{\eta}=\{(i,j) \, | \, 1 \leq i \leq 
\eta_1+\cdots \eta_k< j \leq n \text{ for some
} l \} \, .
\end{equation}
The formal power series is then defined,
for the alphabet $X_n=x_1+\cdots+x_n$, by
\begin{equation}
{\mathcal B}_{\eta}[X_n;t] =
\prod_{(i,j) \in \text{Roots}_{\eta} } \frac{1}{1-tx_i/x_j} \, .
\end{equation}
Given a function $f(x_1,\cdots,x_n)$,
the action of a permutation $\sigma$ of $S_n$ is defined by
\begin{equation}
\sigma \,  f(x_1, x_2, \dots, x_n)  =
f(x_{\sigma(1)} ,x_{\sigma(2)}, \dots, x_{\sigma(n)}) \, .
\end{equation}
If $A$ is the the antisymmetrizer, $\sum_{\sigma \in S_n} \text{sign}(\sigma)\, 
\sigma$,
and $\delta = (n-1,n-2,\dots,0)$,
the operator $\pi$ is defined as
\begin{equation}
\pi(f)= A(x^\delta f)/A(x^{\delta}) \, .
\end{equation}
Denoting the concatenation of the partitions in $S$ by
$\bar S$, we then have the generating series
\begin{equation}
{\mathbb H}_S[X_n;t]= \pi \bigl( \, x^{\bar S}\, {\mathcal B}_{\eta}[X;t] \,
\bigr) \, .
\end{equation}
Since $\pi$ sends a polynomial in the variables of $X_n$
to a symmetric polynomial, $\mathbb H_S$ is symmetric in $X_n$.
Letting $\mathcal P_n$ be the set of elements
of $\mathbb Z^n$ whose entries are weakly decreasing, we can thus write
\begin{equation}
{\mathbb H}_S[X_n;t]=
\sum_{\omega \in {\mathcal {P}_n} } K_{\omega;S}(t) s_{\omega}[X_n]
 \, .
\end{equation}
Considering only partitions of
$|S|\equiv |\lambda^{(1)}| +|\lambda^{(2)}|+\cdots$,
we finally define the symmetric polynomial
\begin{equation}
\H_S[X_n;t]= \sum_{\mu \vdash |S|} K_{\mu;S}(t) s_{\mu}[X_n] \, .
\end{equation}

Note that if $S'$ denotes the sequence of
partitions obtained by appending the partition
$(0^\ell)$ to a sequence of partitions $S$, then $|S'|=|S|$.
In this case, we can show that 
\begin{equation}
\H_{S'}[X_{n+\ell}] = \sum_{\mu \vdash |S'|}
K_{\mu;S'}(t) \, s_{\mu}[X_{n+\ell}]
 = \sum_{\mu \vdash |S|} K_{\mu;S}(t) \, s_{\mu}[X_{n+\ell}] \, ,
\end{equation}
that is, $\H_{S'}[X_{n+\ell}]$ is equal to $\H_{S}[X_{n}\to X_{n+\ell}]$.
The number of variables in $\H_S[X_n]$ is thus irrelevant
(as long as it is large enough) and therefore, at times, we work with the infinite
alphabet $X$.

\subsection{Properties of the generalized Kostka polynomials}
\label{A2}
It is known that the coefficients $K_{\mu;S}(t)$ obey a Morris-type recurrence.
Starting with a one-element sequence $S=(\lambda^{(1)})$, we have
\begin{equation}
\label{morris1}
K_{\mu;S}(t)=
\begin{cases}
1 & \text{ if } \mu=\lambda^{(1)} \\
0 & \text{ otherwise }
\end{cases} \, .
\end{equation}
Then for an arbitrary sequence
$S= (\lambda^{(1)},\lambda^{(2)},\dots)$,
letting $m=\ell(\lambda^{(1)})$ and
$\hat S= (\lambda^{(2)},\lambda^{(3)},\dots)$,
$K_{\mu;S}(t)$ satisfies the recurrence
\begin{equation}
\label{morris}
K_{\mu;S}(t) = \sum_{w \in S_n/(S_m \times S_{n-m})} (-1)^{\ell(w)}
t^{|\alpha(w)|-|\lambda^{(1)}|} \sum_{\nu \vdash |\hat S|}
c^{\nu}_{\alpha(w)/\lambda^{(1)},\beta(w)} K_{\nu;\hat S}(t) \, ,
\end{equation}
where $w$ runs over minimal length coset representatives, and
where $\alpha(w)$ and $\beta(w)$ are the first $m$ and last
$n-m$ parts of the weight $w^{-1}(\mu+\delta)-\delta$,
respectively. Note that $\beta(w)$ is always a partition, and that
the $w$-th summand is understood
to be zero unless all the parts of $\alpha(w)$ are nonnegative and
$\alpha(w) \supseteq \lambda^{(1)}$.
Further, $c^{\nu}_{\alpha(w)/\lambda^{(1)},\beta(w)}$ is
the Littlewood-Richardson coefficient
\begin{equation}
c^{\nu}_{\alpha(w)/\lambda^{(1)},\beta(w)} = \langle s_{\alpha(w)/\lambda^{(1)}}[X]
s_{\beta(w)}[X], s_{\nu}[X] \rangle \, .
\end{equation}
Note that any element $w$ of $S_n/(S_m \times S_{n-m})$
is of the form
\begin{equation}
\label{w}
[[i_1,\dots,i_m]]\, \equiv \,
[i_1,\dots,i_m,1,\dots,\hat i_1,\dots,
\hat i_m,\dots,n]
\qquad 1 \leq i_1 < \cdots < i_m \leq n\, ,
\end{equation}
where $1,\dots,\hat i_1,\dots,\hat i_m,\dots,n$
denotes the sequence $1,\dots,n$ with the elements
$i_1,\dots,i_m$ omitted.
Since any permutation $w=[w_1,\dots,w_n]$
satisfies
\begin{equation}
\label{wm}
w^{-1}(\mu+\delta)-\delta=
\bigl(\mu_{w_1}-(w_1-1),\dots,\mu_{w_n}-(w_n-n)\bigr)
\,
\end{equation}
in the case that $w=[[i_1,\dots,i_m]]$, \eqref{w} implies
that
\begin{equation}
\label{wm1}
(\mu_1,\ldots,\mu_m) \supseteq \alpha(w)
\end{equation}
and
\begin{equation}
\label{wm2}
\mu_j=(\beta(w))_{j-m}\quad\text{for all } \quad j>i_m
\, .
\end{equation}

\medskip

We now prove two properties held by
the generalized Kostka polynomials, $K_{\mu;S}(t)$, when
the sequence $S$ is dominant, that is, when $\bar S$ is a partition.

\medskip

\begin{proposition}
\label{propuni}
If $S$ is a dominant sequence
of partitions with $\bar S=\mu$, then $K_{\mu;S}(t)=1$.
\end{proposition}
\noindent {\bf Proof.}\quad
For $S$ with one element, the claim follows
from (\ref{morris1}).
Let $S=(\lambda^{(1)},\lambda^{(2)},\ldots)$, with $\ell(\lambda^{(1)})=m$.
Since ${\bar S}=\mu$, we have
$K_{\hat \mu;\hat S}=1$, where $\hat \mu=(\mu_{m+1},\mu_{m+2},\dots)$,
 by induction on the number of elements in $S$.
Consider the $w$-th summand of \eqref{morris}.
Since the only non-zero terms
occur when $\alpha(w) \supseteq \lambda^{(1)}=(\mu_1,\dots,\mu_m)$,
and $\alpha(w)\subseteq(\mu_1,\ldots,\mu_m)$
by \eqref{wm1}, the only term is when
$\alpha(w)=\lambda^{(1)}$ and
$\beta(w)=\overline{(\hat S)}\equiv\hat\mu$.
Thus,
$K_{\mu;S}(t)=c^{\nu}_{0,\hat \mu}
K_{\hat \mu;\hat S}=K_{\hat \mu;\hat S}=1$,
since $c^{\nu}_{0,\gamma}=0$ unless $\nu=\gamma$.
\hfill $\endprf$

\medskip

We are also able to prove that many of the generalized
Kostka polynomials vanish.
To this end, we need the following two lemmas.
The first lemma concerns the Littlewood-Richardson coefficients.

\medskip
\begin{lemma} \label{lemLittle}
Let $\nu$, $\alpha$, $\lambda$ and $\beta$ be partitions
and let $N=|\alpha|-|\lambda|$. If
$\beta^{N} \equiv (\beta_1+N,\beta_2,\beta_3,\dots)$
then $c^{\nu}_{\alpha/\lambda,\beta}=0$ when $\nu \not \leq \beta^N$.
\end{lemma}
\noindent {\bf{Proof.}} \quad Given $s_{\alpha/\lambda} 
= \sum_{\rho \vdash N} d_{\rho} \, s_{\rho}$, 
we have, by definition,
$c^{\nu}_{\alpha/\lambda,\beta}= \sum_{\rho \vdash N} d_{\rho} \,
c^{\nu}_{\rho,\beta}$.
Therefore, the lemma will hold if we can show that
$c_{\rho,\beta}^\nu=0$ for all $\rho \vdash N$ when $\nu \not \leq \beta^N$.
Since $\rho \vdash N$ implies  $\rho'\geq (N)'$, we have
$(\rho'\cup\beta')'\leq ((N)'\cup\beta')'=\beta^N$.
But Property~\ref{LR} gives that
$c^{\nu}_{\rho,\beta}=0$ unless $\nu \leq \bigl(\rho' \cup \beta'\bigr)'$.
Thus,
$c^{\nu}_{\rho,\beta}=0$ unless $\nu \leq \bigl(\rho' \cup \beta'\bigr)'
\leq \beta^N$, that is $c_{\rho,\beta}^\nu=0$  when $\nu \not \leq \beta^N$.
\hfill $\endprf$

\medskip

\begin{lemma} \label{lemextra}
Let $\mu$ and $\lambda$ be partitions such that $\mu \not \geq \lambda$, and
let $\alpha$ and $\beta$ denote respectively the first $m$ and
last $n-m$ parts of $w^{-1}(\mu+\delta)-\delta$
for $w \in S_n/(S_m \times S_{n-m})$.
If $\alpha \supseteq (\lambda_{1},\ldots,\lambda_m)$
then $\beta^N \not \geq (\lambda_{m+1},\dots,\lambda_n)$,
where $\beta^N$ is as defined in Lemma~\ref{lemLittle}.
\end{lemma}
\noindent{\bf Proof.}\quad
Since $w=[[i_1,\ldots,i_m]]$, 
\eqref{wm1} implies $(\mu_1,\ldots,\mu_m)\supseteq\alpha$.
Thus, given $\alpha\supseteq (\lambda_1,\ldots,\lambda_m)$,
\begin{equation}
\label{L1}
\lambda_r\leq\alpha_r\leq\mu_r\quad
\text{for all}\quad r\leq m\,.
\end{equation}
Now, $\mu\not\geq\lambda$ implies that
there exists some $i$ that is the
smallest integer such that
\begin{equation}
\label{plu}
\mu_1+\cdots+\mu_i < \lambda_1+\cdots+\lambda_i
\, .
\end{equation}
This necessarily implies that $\mu_i<\lambda_i$,
and therefore $i>m$ by \eqref{L1}.
We now prove a result from which the
lemma will then follow.  That is,
beyond the $i^{th}$ entry, $\alpha \cup \beta$
is identical to $\mu$;
\begin{equation} \label{reorder}
\alpha \cup \beta=(-,\dots,-,\beta_{i-m+1},\dots,\beta_{n-m})
=(-,\dots,-,\mu_{i+1},\dots,\mu_n) \, .
\end{equation}
The lemma follows since $|\alpha|+|\beta|=|\alpha\cup\beta|=|\mu|$ 
implies by (\ref{reorder}) that
$|\alpha|+\beta_1+\cdots+\beta_{i-m}
= (\alpha\cup\beta)_1+\cdots+ (\alpha\cup\beta)_i
= \mu_1+\cdots+\mu_i$.
Which is to say
that $|\alpha| +\beta_1+\cdots+\beta_{i-m}
< \lambda_1+\cdots+\lambda_m+
\lambda_{m+1}+\cdots+\lambda_i$, by \eqref{plu}.
With $N=|\alpha|-(\lambda_1+\cdots+\lambda_m)$, we then
have $N +\beta_1+\cdots+\beta_{i-m}
<\lambda_{m+1}+\cdots+\lambda_i$, which gives
$\beta^N \not \geq (\lambda_{m+1},\dots,\lambda_n)$.

\noindent {\it Proof of (\ref{reorder}):}  
We first show that beyond the $i^{th}$ entry, $w^{-1}(\mu+\delta)-\delta$
and $\mu$ are identical;
\begin{equation} \label{noreorder}
w^{-1}(\mu+\delta)-\delta=(-,\dots,-,\beta_{i-m+1},\dots,\beta_{n-m})
=(-,\dots,-,\mu_{i+1},\dots,\mu_n) \, .
\end{equation}
With $\mu_i<\lambda_i$
and $i>m$, we have $\mu_i<\lambda_i\leq\lambda_m$
and thus $\mu_i<\lambda_m\leq \alpha_m$
by \eqref{L1}. Further, using \eqref{wm}, we have that
$\alpha_m=\mu_{i_m}-(i_m-m) \leq \mu_{i_m}$ since $i_m \geq m$.
Thus, $\mu_i< \alpha_m \leq \mu_{i_m}$, which gives $i_m<i$.
Hence, since $\beta_{j-m}=\mu_j$ for all $j>i_m$ by \eqref{wm2},
this also holds for $j > i$, that is, 
$(\beta_{i-m+1},\dots,\beta_{n-m})=(\mu_{i+1},\ldots,\mu_n)$.
Now, given \eqref{noreorder}, we will show (\ref{reorder}).  
Note that $\alpha \cup \beta$ is the rearrangement of $w^{-1}(\mu+\delta)-\delta$.
Since $\beta$ is a partition, the rearrangement will not concern
the entries of $\beta$. 
$\mu_i<\alpha_m$ implies that $\beta_{i-m+1}=\mu_{i+1}\leq\mu_i<\alpha_m$
(i.e.  $\beta_{i-m+1}$
is smaller than the smallest element of $\alpha$),
and we have
$\alpha\cup\beta=(-,\dots,-,\beta_{i-m+1},\dots,\beta_{n-m})
=(-,\dots,-,\mu_{i+1},\dots,\mu_n)$.  
 \hfill$\endprf$

\medskip

We can now prove the following statement concerning
the generalized Kostka polynomials.

\medskip

\begin{proposition}
\label{propsmaller}
If $S$ is  dominant with $\bar S=\lambda$,
then $K_{\mu;S}(t)=0$ for  $\mu \not \geq \lambda$.
\end{proposition}
\noindent {\bf{Proof.}} \quad
If $S=(\lambda^{(1)})$ has only one element,
then (\ref{morris1}) gives that $K_{\mu;S}(t)=0$
when $\mu \not \geq \lambda^{(1)}$.
Assume by induction on the number of parts of $S$
that $K_{\nu;\hat S}(t)=0$ for all $\nu \not \geq \overline{\hat S}$.
We now show that the right hand side of \eqref{morris} is
zero for $\mu \not \geq {\bar S}=\lambda$.  Since the $w$-th summand is non-zero only
when $\alpha(w)\supseteq(\lambda^{(1)})$,
the conditions of Lemma~\ref{lemextra} are satisfied.
Therefore, $\beta(w)^N$ is
such that $\beta(w)^N \not \geq \overline{(\hat S)}$.
Since Lemma~\ref{lemLittle} gives that
$c^{\nu}_{\alpha(w)/\lambda^{(1)},\beta(w)}=0$
when $\nu \not\leq \beta(w)^N$,
the only non-zero terms in the sum
occur when $\nu\leq \beta(w)^N$, with
$\overline{\hat S}\not\leq\beta(w)^N$.
If we assumed $\nu \geq \overline{\hat S}$,
then $\overline{\hat S}\leq \nu \leq\beta(w)^N$
would imply $\overline{\hat S} \leq\beta(w)^N$.  By
contradiction, we have that the non-zero terms occur when 
$\nu \not \geq\overline{\hat S}$.
Our induction hypothesis then proves all terms are zero.
\hfill$\square$

\newpage
\section{Tables} \label{tabs}
In the tables below, we have not included the 
cases when $k \geq |\lambda|$, which, from Property~\ref{kgrandt}, 
simply correspond to the trivial cases $s_{\lambda}[X;t]=s_{\lambda}[X]$. 

\subsection{$k$-Schur functions in terms of Schur functions} \label{tabkschur}

\;\quad\;
\medskip

\begin{center}

\begin{tabular}{|c||c|c|c|}
\hline
$k=2$ & $1^3$ & 21 & 3\\
\hline\hline
$1^3$& 1 & $t$ & \\ \hline
21&   & 1& $t$\\ \hline
\end{tabular}
\qquad
\begin{tabular}{|c ||c|c|c|c|c| }
\hline
$k=2$ & $1^4$ &$21^2$ & $2^2$ & 31 & 4 \\ 
\hline
\hline
$1^4$ & 1 & $t$ & $t^2$ & & \\ \hline
$21^2$ &  & 1 &  & $t$ & \\ \hline
$2^2$ &  &  & $1$ & $t$ & $t^2$ \\ \hline
\end{tabular}

\begin{tabular}
{|c ||c|c|c|c|c|c|c| }
\hline
$k=2$ & 
$1^5$ & $21^3$ & $2^21$ & $31^2$ & 32 & 41 & 5
\\ 
\hline
\hline
{$1^5$ }  
& 1 & $t+t^2$ & $t^2+t^3$  & $t^3$  & $t^4$  &  &
\\
\hline
{$21^3$}
& & 1 & $t$ & $t+t^2$  & $t^2$ & $t^3$  &
\\
\hline
{$2^21$}
&  &  & $1$  & $t$  & $t+t^2$  & $t^2+t^3$  & $t^4$
\\ 
\hline
\end{tabular}
\qquad
\begin{tabular}
{|c||c|c|c|c|c| }
\hline
{$k=3$} & 
{$1^4$} & {$21^2$} &{$2^2$} & {$31$} & {$4$}
 \\ 
\hline
\hline
{$1^4$} & {1} & $t$ & & &
\\
\hline
{$21^2$} & & {1} & &$t$ &
\\
\hline
{$2^2$} & & & {1} & &
\\
\hline
{$31$} & & & & {1} & {$t$}
\\
\hline
\end{tabular}

\begin{tabular}
{|c||c|c|c|c|c|c|c| }
\hline
{$k=3$} & 
{$1^5$} & {$21^3$} &{$2^21$} & {$31^2$} & {$32$} & {41} & 5
 \\ 
\hline
\hline
{$1^5$} & {1} & $t$ & $t^2$ & & & &
\\
\hline
{$21^3$} & & {1} & & $t$  & & &
\\
\hline
{$2^21$} & & & {1} & & $t$ & &
\\
\hline
{$31^2$} & & & &  {1} & & {$t$}  &
\\
\hline
{$32$} & & & & & {1} & $t$ & $t^2$
\\
\hline
\end{tabular}

\begin{tabular}
{|c||c|c|c|c|c|c|c|c|c|c|c| }
\hline
{$k=3$} & 
{$1^6$} & {$21^4$} &{$2^21^2$} & {$2^3$} & {$31^3$} & {$321$} & $3^2$ & $41^2$ &$42$ &$51$&$6$
 \\ 
\hline
\hline
{$1^6$} & {1} & $t$ & $t^2$ & $t^3$ & & & & & & &
\\
\hline
{$21^4$} & & {1} &$t$  & &$t$&$t^2$ &  & & & &
\\
\hline
{$2^21^2$} & & & {1} &  & &$t$  &$t^2$ & & & &
\\
\hline
{$2^3$} & & & & {1} &  & $t$& & &$t^2$ & &
\\
\hline
{$31^3$} & & & & & 1 &  &  & $t$& & &
\\
\hline
{$321$} & & & & &  &1 &  & $t$&$t$ &$t^2$ &
\\
\hline
{$3^3$} & & & & & &  & 1 &  &$t$ &$t^2$ &$t^3$
\\
\hline
\end{tabular}

\begin{tabular}
{|c||c|c|c|c|c|c|c| }
\hline
{$k=4$} & 
{$1^5$} & {$21^3$} &{$2^21$} & {$31^2$} & {$32$} & {41} & 5
 \\ 
\hline
\hline
{$1^5$} & {1} & $t$ & & & & &
\\
\hline
{$21^3$} & & {1} &  &$t$ & & &
\\
\hline
{$2^21$} & & & 1 &  & & &
\\
\hline
{$31^2$} & & & & 1 &  & $t$ & 
\\
\hline
{$32$} & & & & & 1 &  & 
\\
\hline
{$41$} & & & & & & 1 & $t$
\\
\hline
\end{tabular}

\begin{tabular}
{|c||c|c|c|c|c|c|c|c|c|c|c| }
\hline
{$k=4$} & 
{$1^6$} & {$21^4$} &{$2^21^2$} & {$2^3$} & {$31^3$} & {$321$} & $3^2$ & $41^2$ &$42$ &$51$&$6$
 \\ 
\hline
\hline
{$1^6$} & {1} & $t$ & $t^2$ & & & & & & & &
\\
\hline
{$21^4$} & & {1} &  & &$t$& &  & & & &
\\
\hline
{$2^21^2$} & & & {1} &  & &$t$  & & & & &
\\
\hline
{$2^3$} & & & & {1} &  & & & & & &
\\
\hline
{$31^3$} & & & & & 1 &  &  & $t$& & &
\\
\hline
{$321$} & & & & &  &1 &  & & $t$ & &
\\
\hline
{$3^3$} & & & & & &  & 1&  & & &
\\
\hline
{$41^2$} & & & & & &  &  & 1& &$t$&
\\
\hline
{$42$} & & & & & &  & &  &$1$ &$t$ &$t^2$
\\
\hline
\end{tabular}

\end{center}

\subsection{Macdonald polynomials in terms of $k$-Schur functions} \label{tabmac}

\begin{center}
\quad\;\quad

\begin{tabular}{|c||c|c|}
\hline
$k=2$ & $1^3$ & 21 \\
\hline\hline
$1^3$& 1 & $t^2$\\ \hline
21& $q$  & 1\\ \hline
\end{tabular}
\qquad
\begin{tabular}{|c||c|c|c|}
\hline
$k=2$ & $1^4$ & $21^2$ & $2^2$ \\
\hline\hline
$1^4$& 1 & $t^2+t^3$ & $t^4$\\ \hline
$21^2$& $q$ & $1+qt^2$ & $t$\\ \hline
$2^2$& $q^2$ & $q+qt$  & 1\\ \hline
\end{tabular}
\qquad
\begin{tabular}{|c||c|c|c|}
\hline
$k=2$ & $1^5$ & $21^3$ & $2^21$ \\
\hline\hline
$1^5$& 1 & $t^3+t^4$ & $t^6$\\ \hline
$21^3$& $q$ & $1+qt^3$ & $t^2$\\ \hline
$2^21$& $q^2$  & $q+qt$ & 1\\ \hline
\end{tabular}

\begin{tabular}{|c||c|c|c|c|}
\hline
$k=2$ & $1^6$ & $21^4$ & $2^21^2$ & $2^3$ \\
\hline\hline
$1^6$& 1 & $t^3+t^4+t^5$ & $t^6+t^7+t^8$ & $t^9$\\ \hline
$21^4$& $q$ & $1+qt^3+qt^4$ & $t^2+t^3+qt^6$ & $t^4$\\ \hline
$2^21^2$& $q^2$  & $q+qt+q^2t^3$ & $1+qt^2+qt^3$ & $t$\\ \hline
$2^3$& $q^3$ & $q^2+q^2t+q^2t^2$  & $q+qt+qt^2$ & 1\\ \hline
\end{tabular}
\qquad
\begin{tabular}{|c||c|c|c|c|}
\hline
$k=3$ & $1^4$ & $21^2$ & $2^2$ & $31$ \\
\hline\hline
$1^4$& 1 & $t^2+t^3$ & $t^2+t^4$ & $t^5$\\ \hline
$21^2$& $q$ & $1+qt^2$ & $t+qt^2$ & $t^2$\\ \hline
$2^2$& $q^2$ & $q+qt$ & $1+q^2t^2$ & $t$\\ \hline
$31$& $q^3$ & $q+q^2$  & $q+q^2t$ & 1\\ \hline
\end{tabular}

\begin{tabular}{|c||c|c|c|c|c|}
\hline
$k=3$ & $1^5$ & $21^3$ & $2^21$ & $31^2$ & $32$ \\
\hline\hline
$1^5$& 1 & $t^2+t^3+t^4$ & $t^3+t^4+t^5+t^6$ & $t^5+t^6+t^7$ & $t^8$\\ \hline
$21^3$& $q$ & $1+qt^2+qt^3$ & $t+t^2+qt^3+qt^4$ & $t^2+t^3+qt^5$ & $t^4$\\ \hline
$2^21$& $q^2$ & $q+qt+q^2t^2$ & $1+qt+qt^2+q^2t^3$ & $t+qt^2+qt^3$ & $t^2$\\ \hline
$31^2$& $q^3$ & $q+q^2+q^3t^2$  & $q+qt+q^2t+q^2t^2$ & $1+qt^2+q^2t^2$ & $t$\\ \hline
$32$ &$q^4$ &$q^2+q^3+q^3t$ & $q+q^2+q^2t+q^3t$ & $q+qt+q^2t$ & 1 \\ \hline
\end{tabular}

{\small
\begin{tabular}{|c||c|c|c|c|c|c|}
\hline
$k=4$ & $1^5$ & $21^3$ & $2^21$ & $31^2$ & $32$ & 41 \\
\hline\hline
$1^5$ & 1 & $t^2+t^3+t^4$ & $t^2+t^3+t^4+t^5+t^6$ & $t^5+t^6+t^7$ & $t^4+t^5+t^6
+t^7+t^8$ & $t^9$\\ \hline
$21^3$ & $q$ & $1+qt^2+qt^3$ & $t+t^2+qt^2+qt^3+qt^4$ & $t^2+t^3+qt^5$ & 
$t^2+t^3+t^4+qt^4+qt^5$ & $t^5$\\ \hline
$2^21$ & $q^2$ & $q+qt+q^2t^2$ & $1+qt+qt^2+q^2t^2+q^2t^3$ & $t+qt^2+qt^3$ & 
$t+t^2+qt^2+qt^3+q^2t^4$ & $t^3$\\ \hline
$31^2$ & $q^3$ & $q+q^2+q^3t^2$ & $q+qt+q^2t+q^2t^2+q^3t^2$ &
$ 1+qt^2+q^2t^2$ & $t+qt+qt^2+q^2t^2+q^2t^3$ & $t^2$\\ \hline
$32$ & $q^4$ & $q^2+q^3+q^3t$ & $q+q^2+q^2t+q^3t+q^4t^2$ & $q+qt+q^2t$ & 
$1+qt+q^2t+q^2t^2+q^3t^2$ & $t$\\ \hline
$41$ & $q^6$ & $q^3+q^4+q^5$  & $q^2+q^3+q^4+q^4t+q^5t$ & $q+q^2+q^3$  & 
$q+q^2+q^2t+q^3t+q^4t$ & 1\\ \hline
\end{tabular}}
\end{center}

\end{document}